\documentclass[12pt]{article}
\usepackage{amsmath,amssymb}
\usepackage{graphicx}
\newcommand{\eqdef}{\stackrel{\text{def}}{=}}
\newcommand{\eqdefrm}{\stackrel{\text{\rm def}}{=}}
\newcommand{\n}{\nonumber\\}
\newcommand{\bm}{\boldsymbol}
\newcommand{\ignore}[1]{}
\numberwithin{equation}{section}
\newcommand{\Romannumeral}[1]{\uppercase\expandafter{\romannumeral#1}}
\newcommand{\I}{\text{\Romannumeral{1}}}
\newcommand{\II}{\text{\Romannumeral{2}}}

\newtheorem{theo}{\bf Theorem}[section]

\newcommand{\laprod}[2]{\prod_{#1}^{\stackrel{#2}{\longleftarrow}}}
\newcommand{\cFt}{\tilde{\mathcal{F}}}
\newcommand{\cBt}{\tilde{\mathcal{B}}}
\newcommand{\cX}{\mathcal{X}}
\newcommand{\qbinom}[2]{\genfrac{[}{]}{0pt}{}{\,#1\,}{#2}}

\allowdisplaybreaks[4]

\begin{document}

\baselineskip=20pt

\newcommand{\preprint}{
\vspace*{-20mm}
   \begin{flushright}\normalsize \sf
    DPSU-22-1\\
  \end{flushright}}
\newcommand{\Title}[1]{{\baselineskip=26pt
  \begin{center} \Large \bf #1 \\ \ \\ \end{center}}}
\newcommand{\Author}{\begin{center}
  \large \bf Satoru Odake$^{\,a}$ and Ryu Sasaki$^{\,b}$ \end{center}}
\newcommand{\Address}{\begin{center}
     $^a$ Faculty of Science, Shinshu University,
     Matsumoto 390-8621, Japan\\
     $^b$ Department of Physics, Tokyo University of Science,
     Noda 278-8510, Japan
   \end{center}}
\newcommand{\Accepted}[1]{\begin{center}
  {\large \sf #1}\\ \vspace{1mm}{\small \sf Accepted for Publication}
  \end{center}}

\preprint
\thispagestyle{empty}

\Title{``Diophantine'' and Factorisation Properties of Finite Orthogonal
Polynomials in the Askey Scheme}

\Author

\Address
\vspace{1cm}

\begin{abstract}
A new interpretation and applications of the ``Diophantine'' and factorisation
properties of {\em finite} orthogonal polynomials in the Askey scheme are
explored. The corresponding twelve polynomials are the ($q$-)Racah,
(dual, $q$-)Hahn, Krawtchouk and five types of $q$-Krawtchouk.
These ($q$-)hypergeometric polynomials, defined only for the degrees of
$0,1,\ldots,N$, constitute the main part of the eigenvectors of
$N+1$-dimensional tri-diagonal real symmetric matrices, which correspond to
the difference equations governing the polynomials. The {\em monic} versions
of these polynomials all exhibit the ``Diophantine'' and factorisation
properties at higher degrees than $N$. This simply means that these higher
degree polynomials are zero-norm ``eigenvectors'' of the $N+1$-dimensional
tri-diagonal real symmetric matrices. A new type of multi-indexed orthogonal
polynomials belonging to these twelve polynomials could be introduced by
using the higher degree polynomials as the seed solutions of the multiple
Darboux transformations for the corresponding matrix eigenvalue problems.
The shape-invariance properties of the simplest type of the multi-indexed
polynomials are demonstrated. The explicit transformation formulas are
presented.
\end{abstract}

\section{Introduction}
\label{sec:intro}

It is noted for some time that some of the best known orthogonal polynomials,
the Laguerre, Jacobi, Wilson and Askey-Wilson, etc exhibit strange phenomena,
called ``Diophantine'' and factorisation properties \cite{cal,chen-ismail},
when some parameters are tuned at values which invalidate the orthogonality.
In this paper we show that these properties are universally shared, without
tuning parameters, by a not so well known group of orthogonal polynomials,
the {\em monic} versions of finite polynomials in the Askey scheme
\cite{askey,ismail,kls,gasper}.

A simple interpretation or explanation is that the monic versions of these
polynomials satisfy the same $(N+1)\times(N+1)$ matrix eigenvalue problem
\eqref{highmateq} which govern the polynomials \eqref{difeq}--\eqref{comp}.
Therefore at higher $n=N+1+m$ ($m\in\mathbb{Z}_{\ge0}$) degrees, the
polynomials are zero-norm solutions, which means the ``Diophantine'' and
factorisation properties.

This paper is organised as follows.
Section two start with some typical examples of the ``Diophantine'' and
factorisation properties together with the background description. The finite
orthogonal polynomials in the Askey scheme, the main subject, are briefly
introduced as the solutions of certain tri-diagonal real symmetric
$(N+1)\times(N+1)$ matrix eigenvalue problems \cite{os12} in \S\,\ref{sec:fin}.
The explicit data of the twelve polynomials are displayed according to the
five families of the sinusoidal coordinates in \S\,\ref{sec:data}.
In \S\,\ref{sec:zero} we present the main {\bf Theorem\,\ref{theo1}} stating
that the ``Diophantine'' and factorisation is the consequence of the zero
norm nature of the higher degree $N+1+m$ ($m\in\mathbb{Z}_{\ge0}$) monic
polynomials. Starting with the most generic $q$-Racah polynomial, the
explicit expressions of the ``Diophantine'' properties and factorisation are
displayed in \S\,\ref{sec:Kra1}--\S\,\ref{sec:dqK1}. The general setting of the
multi-indexed orthogonal polynomials generated by multiple Darboux
transformations by using the zero-norm solutions is outlined in section three.
The very special cases of the new multi-indexed polynomials constructed by
using $M$ contiguous lowest degree zero-norm solutions are detailed in section
four. In these cases, the multi-indexed polynomials take the same form as the
original with $x$ and some parameters shifted, displaying the so-called shape
invariance. The general transformation rules for the five families are
listed in {\bf Theorem\,\ref{theo:sim}} together with the explicit forms of the
transformation rules of the polynomials in {\bf Theorem\,\ref{theo:42}}.
Corresponding to the new type of shape invariance, the forward $x$-shift
operators are introduced. It is shown in {\bf Theorem\,\ref{theo:43}} that
the multiple applications of the forward $x$-shift operators reproduce the
results of {\bf Theorem\,\ref{theo:42}}. The explicit forms of the
corresponding backward $x$-shift operators are listed.
The final section is for a summary and some comments.

\section{``Diophantine'' Properties and Factorisation}
\label{sec:dio}

Many interesting examples of the ``Diophantine'' properties and factorisation
of various orthogonal polynomials belonging to the Askey scheme are reported
and explained by Calogero and collaborators \cite{cal} (to be cited as \I)
and Ismail and a coauthor \cite{chen-ismail} (to be cited as \II).
Here we list some typical examples, which are stated for the {\em monic}
versions of the named polynomials: (typos in (\II.3.10) and (\II.4.4) are
corrected)
\begin{align*}
  &\text{Racah with }\alpha=-n,\\
  &\qquad\quad p_n(x;-n,\beta,\gamma,\delta)
  =\prod_{k=0}^{n-1}\bigl(x-k(k+\gamma+\delta+1)\bigr),
  \tag{\I.3.29a}\\
  &\text{Jacobi with }\alpha=-n,\\
  &\qquad\quad p_n(x;-n,\beta)=(x-1)^n,
  \tag{\I.3.153}\\
  &\text{Laguerre with }\alpha=-n,\\
  &\qquad\quad p_n(x;-n)=x^n,
  \tag{\I.3.162}\\
  &\text{Wilson with }t_3=1-m-t_4,\ \ 0\leq m\leq n,\\
  &\quad\ \ p_n(x;t_1,t_2,1-m-t_4,t_4)
  =\prod_{k=0}^{m-1}\bigl(x+(t_4+k)^2\bigr)\cdot
  p_{n-m}(x;t_1,t_2,1-t_4,t_4+m),
  \tag{\II.3.10}\\
  &\text{Askey-Wilson with }t_3=q^{1-m}/t_4,\ \ x\eqdef\cos\theta,
  \ \ 0\leq m\leq n,\\
  &\qquad\quad p_n(x;t_1,t_2,q^{1-m}/t_4,t_4)
  =(t_4e^{i\theta}\,;q)_m(t_4e^{-i\theta}\,;q)_m\,
  p_{n-m}(x;t_1,t_2,q/t_4,q^m t_4)\\
  &\phantom{\qquad\quad p_n(x;t_1,t_2,q^{1-m}/t_4,t_4)=}
  \times(-1)^mt_4^{-m}q^{-\frac12m(m-1)}2^{-m}.
  \tag{\II.4.4}
%
%
%
\end{align*}
It should be stressed that in all these examples, some parameters are tuned
to a degree number $n$, in a rather ad-hoc manner.
Therefore the ``Diophantine'' and factorisation properties do not belong to
the polynomials in general but only to the particular degree polynomial to
which the parameters are tuned. Moreover, with those parameter assignments,
the polynomials are no longer orthogonal with each other.
This situation makes it difficult to find satisfactory interpretations of the
``Diophantine'' and factorisation properties and tends to give a wrong
impression that they are of haphazard or unsystematic origin and having
rather peripheral importance.

In this paper we present a different perspective and show that for a certain
group of orthogonal polynomials in the Askey scheme, the ``Diophantine'' and
factorisation properties are essential and inherent within the proper
parameter ranges in which the orthogonality holds. Let us first introduce
the general features of the polynomials belonging to this group.

\subsection{Finite orthogonal polynomials of a discrete variable}
\label{sec:fin}

The polynomials $\{\check{P}_n(x;N,\bm{\lambda})\}$ in the group are also
called finite orthogonal polynomials of a discrete variable \cite{nikiforov}.
They are defined on a finite integer lattice $\cX=\{0,1,\ldots,N\}$ satisfying 
\begin{equation} 
  \sum_{x\in\cX}w(x;N,\bm{\lambda})\check{P}_m(x;N,\bm{\lambda})
  \check{P}_n(x;N,\bm{\lambda})
  =\frac{\delta_{m\,n}}{d_n(N,\bm{\lambda})^2}\quad(m,n\in\cX),
  \label{Pnnorm}
  \end{equation}
with a positive weight function $w(x;N,\bm{\lambda})>0$. Here $\bm{\lambda}$
stands for the set of parameters other than the lattice size $N$.
The group consists of five families according to the type of the
{\em sinusoidal coordinate} $\eta(x;\bm{\lambda})$,
\begin{equation}
  \check{P}_n(x;N,\bm{\lambda})
  =P_n\bigl(\eta(x;\bm{\lambda});N,\bm{\lambda}\bigr),
\end{equation}
in which $P_n(\eta(x;\bm{\lambda});N,\bm{\lambda})$ is a degree $n$ polynomial
in $\eta(x;\bm{\lambda})$. There are five types of the sinusoidal coordinates
for the finite polynomials in the Askey scheme \cite{os12},
\begin{align*}
  &\text{(i)}:\ \eta(x)=x,\quad
  \text{(ii)}:\ \eta(x)=x(x+d),\quad
  \text{(iii)}:\ \eta(x)=1-q^x,\\
  &\text{(iv)}:\ \eta(x)=q^{-x}-1,\quad
  \text{(v)}:\ \eta(x)=(q^{-x}-1)(1-dq^x),\\
  &\eta(0)=0,\quad\eta(x)>0\ \ (x\in\cX\backslash\{0\}).
\end{align*}
In this paper we adopt the normalisation convention
\begin{equation}
  \check{P}_n(0;N,\bm{\lambda})=P_n(0;N,\bm{\lambda})=1,
\end{equation}
except for the monic polynomials to be introduced shortly.

As the other orthogonal polynomials in the Askey scheme
\cite{askey,ismail, kls,gasper}, the polynomials in this group satisfy
second order difference equation
\begin{align}
  &\quad B(x;N,\bm{\lambda})\bigl(\check{P}_n(x;N,\bm{\lambda})
  -\check{P}_n(x+1;N,\bm{\lambda})\bigr)
  +D(x;N,\bm{\lambda})\bigl(\check{P}_n(x;N,\bm{\lambda})
  -\check{P}_n(x-1;N,\bm{\lambda})\bigr)\n
  &=\mathcal{E}(n;\bm{\lambda})
  \check{P}_n(x;N,\bm{\lambda})\quad(x,n\in\cX),
  \label{difeq}
\end{align}
on top of the three term recurrence relations.
Here the coefficients $B(x;N,\bm{\lambda})$ and $D(x;N,\bm{\lambda})$ are
positive in $\cX$ \cite{os12} and vanishes only at the boundary of $\cX$,
\begin{alignat*}{2}
  B(x;N,\bm{\lambda})&>0\ \ (x\in\cX\backslash\{N\}),
  &\quad D(x;N,\bm{\lambda})&>0\ \ (x\in\cX\backslash\{0\}),\\
  B(N;N,\bm{\lambda})&=0,&\quad D(0;N,\bm{\lambda})&=0.
\end{alignat*}
As shown explicitly in \cite{os12}, the difference equation \eqref{difeq} can
be rewritten as a matrix eigenvalue equation in terms of an
$(N+1)\times(N+1)$ tri-diagonal matrix
$\bigl(\widetilde{\mathcal{H}}(N,\bm{\lambda})_{x\,y}\bigr)_{x,y\in\cX}$,
\begin{align}
  &\widetilde{\mathcal{H}}(N,\bm{\lambda})\check{P}_n(x;N,\bm{\lambda})
  =\mathcal{E}(n;\bm{\lambda})\check{P}_n(x;N,\bm{\lambda}),\n
  &\Bigl(\Longleftrightarrow
  \sum_{y\in\cX}\widetilde{\mathcal{H}}(N,\bm{\lambda})_{x\,y}
  \check{P}_n(y;N,\bm{\lambda})
  =\mathcal{E}(n;\bm{\lambda})\check{P}_n(x;N,\bm{\lambda})\Bigr),
  \label{hteq}\\
  &\widetilde{\mathcal{H}}(N,\bm{\lambda})_{x\,x+1}=-B(x;N,\bm{\lambda}),\quad
  \widetilde{\mathcal{H}}(N,\bm{\lambda})_{x\,x-1}=-D(x;N,\bm{\lambda}),\n
  &\widetilde{\mathcal{H}}(N,\bm{\lambda})_{x\,x}
  =B(x;N,\bm{\lambda})+D(x;N,\bm{\lambda}),\quad
  \widetilde{\mathcal{H}}(N,\bm{\lambda})_{x\,y}=0\ \ (|x-y|\ge2),
  \label{hteq1}
\end{align}
with the eigenvalue $\mathcal{E}(n;\bm{\lambda})$.
By a similarity transformation in terms of a positive diagonal matrix
$\Phi(x;N,\bm{\lambda})$, the matrix $\widetilde{\mathcal{H}}(N,\bm{\lambda})$ 
is related to a {\em real symmetric tri-diagonal matrix}
$\mathcal{H}(N,\bm{\lambda})$,
\begin{align}
  &\mathcal{H}\eqdef\Phi\widetilde{\mathcal{H}}\Phi^{-1}\Leftrightarrow 
  \widetilde{\mathcal{H}}=\Phi^{-1}\mathcal{H}\Phi\ \Leftrightarrow
  \mathcal{H}_{x\,y}=\phi_0(x)\widetilde{\mathcal{H}}_{x\,y}\phi_0(y)^{-1},
  \label{Hdef}\\
  &\Phi_{x\,x}=\phi_0(x),\quad\Phi_{x\,y}=0\ \ (x\neq y).
  \label{Phidef}
\end{align}
In these formulas and hereafter the parameter dependence
of $\widetilde{\mathcal{H}}$, $\mathcal{H}$, $\Phi$, $\check{P}_n$, $B$,
$D$, $\eta$, $\mathcal{E}$ etc is omitted occasionally for simplicity of
presentation.
Here a positive function $\phi_0(x)$ on $\cX$ is introduced by the ratios of
$B(x)$ and $D(x+1)$,
\begin{equation}
  \phi_0(x)\eqdef\sqrt{\prod_{y=0}^{x-1}\frac{B(y)}{D(y+1)}}
  \Leftrightarrow
  \frac{\phi_0(x+1)}{\phi_0(x)}=\frac{\sqrt{B(x)}}{\sqrt{D(x+1)}}
  \ \ (x\in\cX\backslash\{N\}),
  \label{phi0def}
\end{equation}
where $\prod_{j=n}^{n-1}*\eqdef 1$ ($\Rightarrow\phi_0(0)=1$).
The tri-diagonal real symmetric matrix $\mathcal{H}$ \eqref{Hdef}, expressed
explicitly as
\begin{align}
  &\mathcal{H}_{x\,x+1}=-\sqrt{B(x)D(x+1)},\quad
  \mathcal{H}_{x\,x-1}=-\sqrt{B(x-1)D(x)},\n
  &\mathcal{H}_{x\,x}=B(x)+D(x),\quad
  \mathcal{H}_{x\,y}=0\ \ (|x-y|\ge2),
  \label{tdefcomp}
\end{align}
has the following eigenvectors
\begin{align}
  &\phi_n(x)\eqdef\phi_0(x)\check{P}_n(x)\ \ (n\in\cX),\quad
  \bigl(\phi_n(x)\bigr)_{x\in\cX}\in\mathbb{R}^{N+1},
  \label{phindef}\\
  &\mathcal{H}\phi_n(x)=\sum_{y\in\cX}\mathcal{H}_{x\,y}\phi_n(y)=
  \sum_{y\in\cX}\phi_0(x)\widetilde{\mathcal{H}}_{x\,y}\check{P}_n(y)
  =\mathcal{E}(n)\phi_n(x)\ \ (n\in\cX).
  \label{comp}
\end{align}
The norms are all finite as
$w(x)\eqdef\sum_{x\in\cX}\phi_0(x)^2<\infty$ for all $B(x)$ and $D(x)$
in the group.
We arrive at the complete set of eigenvectors $\{\phi_n(x)\}$ $(n\in\cX)$,
since the orthogonality of the eigenvectors \eqref{Pnnorm} is guaranteed by
the simpleness of the eigenvalues of the tri-diagonal matrices.

\subsection{Polynomials data}
\label{sec:data}

Here we present the basic data of the polynomials in the group.
For more details of the polynomials, consult \cite{askey}--\cite{gasper}.
The parametrisations of some polynomials are different from the conventional
ones, see \cite{os12}.
The polynomials are divided into five families according to the forms of the
sinusoidal coordinate $\eta(x)$. The ``Diophantine'' properties are identical
within the same family and the factorisation is universal for all the
polynomials in this group.

\paragraph{(i) Family with $\bm{\eta(x)=x}$}
~\\
Two polynomials belong to this family, the Krawtchouk and Hahn polynomials.
They are polynomials in $x$.

\subsubsection{Krawtchouk (K)}
\label{sec:Kra}

The polynomial depends on one positive parameter $\bm{\lambda}=p$ ($0<p<1$).
\begin{align}
  &B(x;N,p)=p(N-x),\quad D(x;N,p)=(1-p)x,\quad
  \mathcal{E}(n)=n,\quad\eta(x)=x,
  \label{Ketadef}\\[2pt]
  &\check{P}_n(x;N,p)=P_n(x;N,p)
  ={}_2F_1\Bigl(\genfrac{}{}{0pt}{}{-n,\,-x}{-N}\Bigm|p^{-1}\Bigr),\quad
  P_n(x;N,p)=P_x(n;N,p),
  \label{PnK}\\[2pt]
  &\check{P}_n(N-x;N,p)=(-1)^n(p^{-1}-1)^n\check{P}_n(x;N,1-p)
  \ \ (\text{Mirror symmetry}).
  \label{Kpmirror}
\end{align}

\subsubsection{Hahn (H)}
\label{sec:Han}

The Hahn polynomial depends on two positive parameters $\bm{\lambda}=(a,b)$
($a,b>0$).
\begin{align}
  &B(x;\bm{\lambda})=(x+a)(N-x),\quad D(x;\bm{\lambda})= x(b+N-x),
  \label{hahnBD}\\
  &\mathcal{E}(n;\bm{\lambda})=n(n+a+b-1),\quad\eta(x)=x,\n
  &\check{P}_n(x;N,\bm{\lambda})=P_n(x;N,\bm{\lambda})
  ={}_3F_2\Bigl(\genfrac{}{}{0pt}{}{-n,\,n+a+b-1,\,-x}{a,\,-N}\Bigm|1\Bigr),
  \label{PnH}\\
  &\check{P}_n(N-x;N,a,b)=(-1)^n\frac{(b)_n}{(a)_n}\check{P}_n(x;N,b,a)
  \ \ (\text{Mirror symmetry}).
  \label{Hpd}
\end{align}

\paragraph{(ii) Family with $\bm{\eta(x)=x(x+d)}$}
~\\
Two polynomials belong to this family, the Racah and dual Hahn polynomials.
They are polynomials in $x(x+d)$.

\subsubsection{Racah (R)}
\label{sec:Rac}

The Racah polynomial depends on three parameters $\bm{\lambda}=(b,c,d)$
($0<d<b-N$, $0<c<1+d$) on top of the lattice size $N$.
\begin{align}
  &B(x;N,\bm{\lambda})=\frac{(x+b)(x+c)(x+d)(N-x)}{(2x+d)(2x+1+d)},\\
  &D(x;N,\bm{\lambda})=\frac{(b-d-x)(x+d-c)x(x+d+N)}{(2x-1+d)(2x+d)},
  \label{RacBD}\\
  &\mathcal{E}(n;\bm{\lambda})=n(n+\tilde{d}),\quad
  \eta(x;\bm{\lambda})=x(x+d),\quad\tilde{d}\eqdef b+c-d-N-1,\n
  &\check{P}_n(x;N,\bm{\lambda})
  =P_n\bigl(\eta(x;\bm{\lambda});N,\bm{\lambda}\bigr)
  ={}_4F_3\Bigl(\genfrac{}{}{0pt}{}{-n,\,n+\tilde{d},\,-x,\,x+d}
  {b,\,c,\,-N}\Bigm|1\Bigr).
  \label{PnR}
\end{align}

\subsubsection{dual Hahn (dH)}
\label{sec:dHn}

The dual Hahn polynomial depends on two positive parameters
$\bm{\lambda}=(a,b)$ ($a,b>0$), and the parameter $d$ is $d=a+b-1$.
\begin{align}
  &B(x;N,\bm{\lambda})=\frac{(x+a)(x+a+b-1)(N-x)}{(2x-1+a+b)(2x+a+b)},
  \label{dualhahnBD1}\\
  &D(x;N,\bm{\lambda})=\frac{x(x+b-1)(x+a+b+N-1)}{(2x-2+a+b)(2x-1+a+b)},
  \label{dualhahnBD2}\\
  &\mathcal{E}(n)=n,\quad\eta(x;\bm{\lambda})=x(x+a+b-1),\n
  &\check{P}_n(x;N,\bm{\lambda})
  =P_n\bigl(\eta(x;\bm{\lambda});N,\bm{\lambda}\bigr)
  ={}_3F_2\Bigl(\genfrac{}{}{0pt}{}{-n,\,x+a+b-1,\,-x}{a,\,-N}\Bigm|1\Bigr).
  \label{PndH}
\end{align}

\medskip

The remaining three families belong to the $q$-hypergeometric polynomial
category. The parameter $q$, $0<q<1$, dependence is not explicitly displayed. 

\paragraph{(iii) Family with $\bm{\eta(x)=1-q^x}$}


\subsubsection{dual quantum $\bm{q}$-Krawtchouk (dq$\bm{q}$K)}
\label{sec:dqqK}

The only polynomial belonging to this family is the dual quantum
$q$-Krawtchouk polynomial \cite{os12} depending on one positive parameter
$\bm{\lambda}=p>q^{-N}$.
\begin{align}
  &B(x;N,p)=p^{-1}q^{-x-N-1}(1-q^{N-x}),\quad
  D(x;N,p)=(q^{-x}-1)(1-p^{-1}q^{-x}),
  \label{dualqqkrawt1}\\
  &\mathcal{E}(n)=q^{-n}-1,\quad\eta(x)=1-q^x,\n
 &\check{P}_n(x;N,p)=P_n\bigl(\eta(x);N,p\bigr)
  ={}_2\phi_1\Bigl(\genfrac{}{}{0pt}{}{q^{-n},\,q^{-x}}{q^{-N}}
  \Bigm|q\,;pq^{x+1}\Bigr).
  \label{PndqqK}
\end{align}

\paragraph{(iv) Family with $\bm{\eta(x)=q^{-x}-1}$}
~\\
On top of the $q$-Hahn polynomial, three sibling polynomials belong to this
family, the $q$-Krawtchouk, quantum $q$-Krawtchouk and affine $q$-Krawtchouk
polynomials, all depending on one positive parameter $\bm{\lambda}=p>0$,
but the parameter ranges are different, as shown in each entry.

\subsubsection{$\bm{q}$-Hahn ($\bm{q}$H)}
\label{sec:qH}

The $q$-Hahn polynomial depends on two positive parameters
$\bm{\lambda}=(a,b)$ ($0<a,b<1$).
\begin{align}
  &B(x;N,\bm{\lambda})=(1-aq^x)(q^{x-N}-1),\quad
  D(x;N,\bm{\lambda})=aq^{-1}(1-q^x)(q^{x-N}-b),
  \label{qHBD}\\
  &\mathcal{E}(n;\bm{\lambda})=(q^{-n}-1)(1-abq^{n-1}),\quad\eta(x)=q^{-x}-1,\n
  &\check{P}_n(x;N,\bm{\lambda})=P_n\bigl(\eta(x);N,\bm{\lambda}\bigr)
  ={}_3\phi_2\Bigl(\genfrac{}{}{0pt}{}{q^{-n},\,abq^{n-1},\,q^{-x}}
  {a,\,q^{-N}}\Bigm|q\,;q\Bigr).
  \label{PnqH}
\end{align}

\subsubsection{$\bm{q}$-Krawtchouk ($\bm{q}$K)}
\label{sec:qK}

\begin{align}
  &B(x;N,p)=q^{x-N}-1,\quad D(x;N,p)=p(1-q^x),
  \label{qKBD}\\
  &\mathcal{E}(n;p)=(q^{-n}-1)(1+pq^n),\quad\eta(x)=q^{-x}-1,\quad p>0,\n
  &\check{P}_n(x;N,p)=P_n\bigl(\eta(x);N,p\bigr)
  ={}_3\phi_2\Bigl(
  \genfrac{}{}{0pt}{}{q^{-n},\,q^{-x},\,-pq^n}{q^{-N},\,0}\Bigm|q\,;q\Bigr).
  \label{PnqK}
\end{align}

\subsubsection{quantum $\bm{q}$-Krawtchouk (q$\bm{q}$K)}
\label{sec:qqK}

\begin{align}
  &B(x;N,p)=p^{-1}q^x(q^{x-N}-1),\quad D(x;N,p)=(1-q^x)(1-p^{-1}q^{x-N-1}),
  \label{qqKBD}\\
  &\mathcal{E}(n)=1-q^n,\quad\eta(x)=q^{-x}-1,\quad p>q^{-N},\n
  &\check{P}_n(x;N,p)=P_n\bigl(\eta(x);N,p\bigr)
  ={}_2\phi_1\Bigl(
  \genfrac{}{}{0pt}{}{q^{-n},\,q^{-x}}{q^{-N}}\Bigm|q\,;pq^{n+1}\Bigr).
  \label{PnqqK}
\end{align}

\subsubsection{affine $\bm{q}$-Krawtchouk (a$\bm{q}$K)}
\label{sec:aqK}

\begin{align}
  &B(x;N,p)=(q^{x-N}-1)(1-pq^{x+1}),\quad D(x;N,p)=pq^{x-N}(1-q^x),\\
  &\mathcal{E}(n)=q^{-n}-1,\quad\eta(x)=q^{-x}-1,\quad 0<p<q^{-1},\n
  &\check{P}_n(x;N,p)=P_n\bigl(\eta(x);N,p\bigr)
  ={}_3\phi_2\Bigl(
  \genfrac{}{}{0pt}{}{q^{-n},\,q^{-x},\,0}{pq,\,q^{-N}}\Bigm|q\,;q\Bigr).
  \label{PnaqK}
\end{align}

\paragraph{(v) Family with $\bm{\eta(x)=(q^{-x}-1)(1-dq^x)}$}
~\\
Three polynomials belong to this family, the $q$-Racah, dual $q$-Hahn
and dual $q$-Krawtchouk polynomials.

\subsubsection{$\bm{q}$-Racah ($\bm{q}$R)}
\label{sec:qRa}

The $q$-Racah polynomial depends on three parameters $\bm{\lambda}=(b,c,d)$
($0<bq^{-N}<d<1$, $qd<c<1$, $\tilde{d}\eqdef bcd^{-1}q^{-N-1}$).
\begin{align}
  &B(x;N,\bm{\lambda})=\frac{(1-bq^x)(1-cq^x)(1-dq^x)(q^{x-N}-1)}
  {(1-dq^{2x})(1-dq^{2x+1})}\,,\\
  &D(x;N,\bm{\lambda})
  =\tilde{d}\,\frac{(b^{-1}dq^x-1)(1-c^{-1}dq^x)(1-q^x)(1-dq^{x+N})}
  {(1-dq^{2x-1})(1-dq^{2x})},\\
  &\mathcal{E}(n;\bm{\lambda})=(q^{-n}-1)(1-\tilde{d}q^n),\quad
  \eta(x;\bm{\lambda})=(q^{-x}-1)(1-dq^x),\n
  &\check{P}_n(x;N,\bm{\lambda})
  =P_n\bigl(\eta(x;\bm{\lambda});N,\bm{\lambda}\bigr)
  ={}_4\phi_3\Bigl(\genfrac{}{}{0pt}{}{q^{-n},\,\tilde{d}q^n,\,q^{-x},\,dq^x}
  {b,\,c,\,q^{-N}}\Bigm|q\,;q\Bigr).
  \label{PnqR}
\end{align}

\subsubsection{dual $\bm{q}$-Hahn (d$\bm{q}$H)}
\label{sec:dqH}

The dual $q$-Hahn polynomial depends on two parameters $\bm{\lambda}=(a,b)$
($0<a,b<1$), and the parameter $d$ is $d=abq^{-1}$.
\begin{align}
  &B(x;N,\bm{\lambda})
  =\frac{(1-aq^x)(1-abq^{x-1})(q^{x-N}-1)}{(1-abq^{2x-1})(1-abq^{2x})},\\
  &D(x;N,\bm{\lambda})=aq^{x-N-1}
  \frac{(1-q^x)(1-bq^{x-1})(1-abq^{x+N-1})}{(1-abq^{2x-2})(1-abq^{2x-1})},\\
  &\mathcal{E}(n)=q^{-n}-1,\quad\eta(x;\bm{\lambda})=(q^{-x}-1)(1-abq^{x-1}),\n
  &\check{P}_n(x;N,\bm{\lambda})
  =P_n\bigl(\eta(x;\bm{\lambda});N,\bm{\lambda}\bigr)
  ={}_3\phi_2\Bigl(\genfrac{}{}{0pt}{}{q^{-n},\,abq^{x-1},\,q^{-x}}
  {a,\,q^{-N}}\Bigm|q\,;q\Bigr).
  \label{PndqH}
\end{align}

\subsubsection{dual $\bm{q}$-Krawtchouk (d$\bm{q}$K)}
\label{sec:dqK}

Similar to the $q$-Krawtchouk, the dual $q$-Krawtchouk polynomial depends on
one parameter $\bm{\lambda}=p>0$, and the parameter $d$ is $d=-p$.
\begin{align}
  &B(x;N,p)=\frac{(q^{x-N}-1)(1+pq^x)}{(1+pq^{2x})(1+pq^{2x+1})},
  \ \ D(x;N,p)
  =pq^{2x-N-1}\frac{(1-q^x)(1+pq^{x+N})}{(1+pq^{2x-1})(1+pq^{2x})},\\
  &\mathcal{E}(n)=q^{-n}-1,\quad\eta(x;p)=(q^{-x}-1)(1+pq^x),\n
  &\check{P}_n(x;N,p)=P_n\bigl(\eta(x;p);N,p\bigr)
  ={}_3\phi_2\Bigl(
  \genfrac{}{}{0pt}{}{q^{-n},\,q^{-x},\,-pq^x}{q^{-N},\,0}\Bigm|q\,;q\Bigr).
  \label{PndqK}
\end{align}

\subsection{Zero norm $=$ ``Diophantine'' and factorisation}
\label{sec:zero}

The polynomials listed above \S\,\ref{sec:Kra}--\S\,\ref{sec:dqK} are
{\em finite} polynomials, since their degree $n$ is limited in
$\cX\ni n$ due to the presence of the lower indices $-N$ or $q^{-N}$
in the ($q$)-hypergeometric expressions, ${}_{r+1}F_r$, ${}_{r+1}\phi_r$,
of these polynomials,
\eqref{PnK}, \eqref{PnH}, \eqref{PnR}, \eqref{PndH}, \eqref{PndqqK},
\eqref{PnqH},
\eqref{PnqK}, \eqref{PnqqK}, \eqref{PnaqK}, \eqref{PnqR}, \eqref{PndqH}
and \eqref{PndqK}. For non-negative integer $m\ge0$, the polynomial
$\check{P}_{N+1+m}(x;N,\bm{\lambda})$ is ill-defined due the presence of the
factors $(-N)_{N+1+k}=0$ and $(q^{-N}\,;q)_{N+1+k}=0$, $0\leq k\leq m$ in the
denominator of the ($q$)-hypergeometric series expansion.

The situation is drastically changed by the introduction of the monic version
of each polynomial in \S\,\ref{sec:data},
\begin{equation}
  \check{P}^{\text{monic}}_n(x;N,\bm{\lambda})
  =\frac1{c_n(N,\bm{\lambda})}\check{P}_n(x;N,\bm{\lambda}),
  \label{cPnmonic}
\end{equation}
where $c_n(N,\bm{\lambda})$ is the coefficient of the highest degree term
$\eta(x;\bm{\lambda})^n$.
It is well defined for all degrees $n=N+1+m$ ($m\in\mathbb{Z}_{\ge0}$) and it
satisfies the same difference equation as \eqref{difeq}, \eqref{hteq}
\begin{align}
  &\qquad\,
  B(x;N,\bm{\lambda})\bigl(\check{P}^{\text{monic}}_{N+1+m}(x;N,\bm{\lambda})
  -\check{P}^{\text{monic}}_{N+1+m}(x+1;N,\bm{\lambda})\bigr)\n
  &\quad+
  D(x;N,\bm{\lambda})\bigl(\check{P}^{\text{monic}}_{N+1+m}(x;N,\bm{\lambda})
  -\check{P}^{\text{monic}}_{N+1+m}(x-1;N,\bm{\lambda})\bigr)\n
  &=\mathcal{E}(N+1+m;\bm{\lambda})
  \check{P}^{\text{monic}}_{N+1+m}(x;N,\bm{\lambda})\quad(x\in\cX),
  \label{difeq1}\\
  &\Longrightarrow
  \widetilde{\mathcal{H}}(N,\bm{\lambda})
  \check{P}^{\text{monic}}_{N+1+m}(x;N,\bm{\lambda})
  =\mathcal{E}(N+1+m;\bm{\lambda})
  \check{P}^{\text{monic}}_{N+1+m}(x;N,\bm{\lambda}).
  \label{hteq2}
\end{align}
This also means that the corresponding monic vector
\begin{equation}
  \phi^{\text{monic}}_{N+1+m}(x)\eqdef
  \phi_0(x)\check{P}^{\text{monic}}_{N+1+m}(x)\quad(m\in\mathbb{Z}_{\ge0}),
  \label{phindef1}
\end{equation}
satisfies the eigenvalue equation of the tri-diagonal real symmetric matrix
$\mathcal{H}$,
\begin{equation}
  \mathcal{H}\phi^{\text{monic}}_{N+1+m}(x)
  =\mathcal{E}(N+1+m;\bm{\lambda})\phi^{\text{monic}}_{N+1+m}(x)
  \quad(m\in\mathbb{Z}_{\ge0}).
  \label{highmateq}
\end{equation}
If $\check{P}^{\text{monic}}_{N+1+m}(x)$ takes a non-vanishing value at some
point in $\cX$, $\phi^{\text{monic}}_{N+1+m}(x)$ becomes an eigenvector of
$\mathcal{H}$, which contradicts the complete set of eigenvectors
$\{\phi_n(x)\}$ $(n\in\cX)$, \eqref{comp}.
Thus we arrive at
\begin{theo}
\label{theo1}
For all the finite polynomials listed in \S\,\ref{sec:data}, the ``Diophantine''
property and factorisation hold,
\begin{align}
  \check{P}^{\text{\rm monic}}_{N+1+m}(x;N,\bm{\lambda})
  &=\Lambda(x;N,\bm{\lambda})\check{Q}_m(x;N,\bm{\lambda})\quad
  (m\in\mathbb{Z}_{\ge0}),
  \label{Dfac}\\
  \Lambda(x;N,\bm{\lambda})&\eqdefrm
  \prod_{k=0}^N\bigl(\eta(x;\bm{\lambda})-\eta(k;\bm{\lambda})\bigr),
  \label{lamdef}
\end{align}
in which $\check{Q}_m(x;N,\bm{\lambda})$
is a monic degree $m$ polynomial in $\eta(x;\bm{\lambda})$.
\end{theo}

Of course, ``Diophantine'' property and factorisation can be demonstrated for
each polynomial by direct calculation of the ($q$)-hypergeometric series
expansion. Below we show the explicit derivation of the ``Diophantine'' and
factorisation property for the $q$-Racah polynomial, the most generic member
of the group.

Since the coefficient of the highest degree of
$\check{P}_n(x\,;N,\bm{\lambda})$ \eqref{PnqR} is
\begin{equation}
  c_n(N,\bm{\lambda})=\frac{(\tilde{d}q^n\,;q)_n}{(b,c,q^{-N}\,;q)_n},
\end{equation}
the monic $q$-Racah polynomial is
\begin{equation}
  \check{P}^{\text{monic}}_n(x;N,\bm{\lambda})
  =\frac{1}{c_n(N,\bm{\lambda})}\check{P}_n(x;N,\bm{\lambda})
  =\sum_{k=0}^n\frac{(bq^k,cq^k,q^{k-N}\,;q)_{n-k}}
  {(\tilde{d}q^{n+k}\,;q)_{n-k}}\frac{(q^{-n},q^{-x},dq^x\,;q)_k}{(q\,;q)_k}q^k.
  \label{cPmonic}
\end{equation}
At $n=N+1+m$ ($m\in\mathbb{Z}_{\ge0}$), $(q^{k-N}\,;q)_{n-k}$ vanishes for
$k=0,1,\ldots,N$, and the $k$ summation is reduced to $\sum_{k=N+1}^{N+1+m}$.
By changing $k=N+1+l$, we obtain
\begin{equation*}
  \check{P}^{\text{monic}}_{N+1+m}(x;N,\bm{\lambda})
  =\sum_{l=0}^m\frac{(bq^{N+1+l},cq^{N+1+l},q^{l+1}\,;q)_{m-l}}
  {(\tilde{d}q^{2(N+1)+m+l}\,;q)_{m-l}}
  \frac{(q^{-N-1-m},q^{-x},dq^x\,;q)_{N+1+l}}{(q\,;q)_{N+1+l}}q^{N+1+l}.
\end{equation*}
{}From the explicit form of $\eta(x;\bm{\lambda})$
\begin{equation*}
  \eta(x;\bm{\lambda})-\eta(k;\bm{\lambda})=-q^{-k}(1-q^{-x+k})(1-dq^{x+k}),
\end{equation*}
we obtain
\begin{equation}
  \Lambda(x;N,\bm{\lambda})=
  \prod_{k=0}^N\bigl(\eta(x;\bm{\lambda})-\eta(k;\bm{\lambda})\bigr)
  =(-1)^{N+1}q^{-\binom{N+1}{2}}(q^{-x},dq^x\,;q)_{N+1}.
  \label{prodeta}
\end{equation}
By using the basic properties of the $q$-shifted factorial,
\begin{equation}
  (a\,;q)_{n+k}=(a\,;q)_n(aq^n\,;q)_k,\quad
  (a\,;q)_n=(-a)^nq^{\binom{n}{2}}(a^{-1}q^{1-n}\,;q)_n,
  \label{qpoch}
\end{equation}
we obtain
\begin{equation*}
  \frac{(q^{-N-1-m},q^{-x},dq^x\,;q)_{N+1+l}}{(q\,;q)_{N+1+l}}
  =\frac{(q^{-N-1-m},q^{-x},dq^x\,;q)_{N+1}}{(q\,;q)_l}
  \frac{(q^{-m},q^{-x+N+1},dq^{x+N+1}\,;q)_l}{(q^{l+1}\,;q)_{N+1}}
\end{equation*}
and
\begin{align*}
  \frac{(q^{-N-1-m}\,;q)_{N+1}}{(q^{l+1}\,;q)_{N+1}}
  &=\frac{(-q^{-N-1-m})^{N+1}q^{\binom{N+1}{2}}(q^{m+1}\,;q)_{N+1}}
  {(q^{l+1}\,;q)_{N+1}}\\
  &=(-q^{-N-1-m})^{N+1}q^{\binom{N+1}{2}}
  \frac{(q^{N+l+2}\,;q)_{m-l}}{(q^{l+1}\,;q)_{m-l}}.
\end{align*}
These lead to
\begin{align}
  \check{P}^{\text{monic}}_{N+1+m}(x;\bm{\lambda})
  &=\Lambda(x;N,\bm{\lambda})\cdot q^{-(N+1)m}\n
  &\quad\times\sum_{l=0}^m
  \frac{(bq^{N+1+l},cq^{N+1+l},qq^{N+1+l}\,;q)_{m-l}}
  {(\tilde{d}q^{2(N+1)+m+l}\,;q)_{m-l}}
  \frac{(q^{-m},q^{-x+N+1},dq^{x+N+1}\,;q)_l}{(q\,;q)_l}q^l\n
  &=\Lambda(x;N,\bm{\lambda})\cdot q^{-(N+1)m}
  \check{P}^{\text{monic}}_m(x-N-1;-N-2,\bm{\lambda}').
  \label{Pmonic4}
\end{align}
in which the shifted parameters are
\begin{equation}
  \bm{\lambda}'=(bq^{N+1},cq^{N+1},dq^{2(N+1)}).
  \label{la'}
\end{equation}

Below we show the ``Diophantine'' and factorisation property for the other
polynomials in the group.

\paragraph{(i) Family with $\bm{\eta(x)=x}$}
~\\
This family possesses the true Diophantine property as all the zeros of
$\Lambda(x)$ are integers. Miki, Tsujimoto and Vinet showed the Diophantine
and factorisation property of the Krawtchouk polynomial by explicit
calculations in the seminal paper \cite{miki}.

\subsubsection{Krawtchouk (K)}
\label{sec:Kra1}

\begin{align}
  &\check{P}^{\text{monic}}_{N+1+m}(x;N,p)
  =\Lambda(x;N)\check{P}^{\text{monic}}_m(x-N-1;-N-2,p),\quad
  c_n(N,p)=\frac{1}{(-N)_n\,p^n},
  \label{PmoKr}\\
  &\check{P}^{\text{monic}}_m(x-N-1;-N-2,p)
  =\sum_{k=0}^m(N+2+k)_{m-k}\frac{(-m,-x+N+1)_k}{k!}p^{m-k}.
  \label{PmoKr1}
\end{align}

\subsubsection{Hahn (H)}
\label{sec:Han1}

\begin{align}
  &\check{P}^{\text{monic}}_{N+1+m}(x;N,\bm{\lambda})
  =\Lambda(x;N)
  \check{P}^{\text{monic}}_m(x-N-1;-N-2,\bm{\lambda}'),
  \label{PmoHa}\\
  &\qquad\bm{\lambda}'=(a+N+1,b+N+1),\quad
  c_n(N,\bm{\lambda})=\frac{(n+a+b-1)_n}{(a,-N)_n},\\
  &\check{P}^{\text{monic}}_m(x-N-1;-N-2,\bm{\lambda}')
  =\sum_{k=0}^m\frac{(a+N+1+k,N+2+k)_{m-k}}{(m+a+b+2N+1+k)_{m-k}}
  \frac{(-m,-x+N+1)_k}{k!}.
  \label{PmoHa1}
\end{align}

\paragraph{(ii) Family with $\bm{\eta(x)=x(x+d)}$}

\subsubsection{Racah (R)}
\label{sec:Rac1}

\begin{align}
  &\check{P}^{\text{monic}}_{N+1+m}(x;N,\bm{\lambda})
  =\Lambda(x;N,\bm{\lambda})
  \check{P}^{\text{monic}}_m(x-N-1;-N-2,\bm{\lambda}'),
  \label{PmoRa}\\
  &\qquad\bm{\lambda}'=\bigl(b+N+1,c+N+1,d+2(N+1)\bigr),\quad
  c_n(N,\bm{\lambda})=\frac{(\tilde{d}+n)_n}{(b,c,-N)_n},\\
  &\check{P}^{\text{monic}}_m(x-N-1;-N-2,\bm{\lambda}')
  =\sum_{k=0}^m\frac{(b+N+1+k,c+N+1+k,N+2+k)_{m-k}}
  {(\tilde{d}+2N+2+m+k)_{m-k}}\n
  &\phantom{\check{P}^{\text{monic}}_m(x-N-1;-N-2,\bm{\lambda}')=\sum_{k=0}^m}
  \times\frac{(-m,-x+N+1,x+N+1+d)_k}{k!}.
  \label{PmoRa1}
\end{align}

\subsubsection{dual Hahn (dH)}
\label{sec:dHn1}

\begin{align}
  &\check{P}^{\text{monic}}_{N+1+m}(x;N,\bm{\lambda})
  =\Lambda(x;N,\bm{\lambda})
  \check{P}^{\text{monic}}_m(x-N-1;-N-2,\bm{\lambda}'),
  \label{PmodH}\\
  &\qquad\bm{\lambda}'=\bigl(a+N+1,b+N+1\bigr),\quad
  c_n(N,\bm{\lambda})=\frac{1}{(a,-N)_n},\\
  &\check{P}^{\text{monic}}_m(x-N-1;-N-2,\bm{\lambda}')
  =\sum_{k=0}^m(a+N+1+k,N+2+k)_{m-k}\n
  &\phantom{\check{P}^{\text{monic}}_m(x-N-1;-N-2,\bm{\lambda}')=\sum_{k=0}^m}
  \times\frac{(-m,-x+N+1,x+a+b+N)_k}{k!}.
  \label{PmodH1}
\end{align}

\paragraph{(iii) Family with $\bm{\eta(x)=1-q^x}$}

\subsubsection{dual quantum $\bm{q}$-Krawtchouk (dq$\bm{q}$K)}
\label{sec:dqqK1}

\begin{align}
  &\check{P}^{\text{monic}}_{N+1+m}(x;N,p)
  =\Lambda(x;N)\,q^{(N+1)m}\check{P}^{\text{monic}}_m(x-N-1;-N-2,p'),
  \label{PmodqqK}\\
  &\qquad p'=pq^{N+1},\quad
  c_n(N,p)=\frac{p^nq^{-\frac12n(n-1)}}{(q^{-N}\,;q)_n},\\
  &\check{P}^{\text{monic}}_m(x-N-1;-N-2,p')
  =\sum_{k=0}^m(q^{N+2+k}\,;q)_{m-k}
  \frac{(q^{-m},q^{-x+N+1}\,;q)_k}{(q\,;q)_k}\n
  &\phantom{\check{P}^{\text{monic}}_m(x-N-1;-N-2,p')=\sum_{k=0}^m}
  \times p^{k-m}q^{kx+k+\frac12m(m-1)-m(N+1)}.
  \label{PmodqqK1}
\end{align}

\paragraph{(iv) Family with $\bm{\eta(x)=q^{-x}-1}$}

\subsubsection{$\bm{q}$-Hahn ($\bm{q}$H)}
\label{sec:qH1}

\begin{align}
  &\check{P}^{\text{monic}}_{N+1+m}(x;N,\bm{\lambda})
  =\Lambda(x;N)\,q^{-(N+1)m}
  \check{P}^{\text{monic}}_m(x-N-1;-N-2,\bm{\lambda}'),
  \label{PmoqH}\\
  &\qquad\bm{\lambda}'=(aq^{N+1},bq^{N+1}),\quad
  c_n(N,\bm{\lambda})=\frac{(abq^{n-1}\,;q)_n}{(a,q^{-N}\,;q)_n},\\
  &\check{P}^{\text{monic}}_m(x-N-1;-N-2,\bm{\lambda}')
  =\sum_{k=0}^m
  \frac{(aq^{N+1+k},q^{N+2+k}\,;q)_{m-k}}{(abq^{m+2N+1+k}\,;q)_{m-k}}
  \frac{(q^{-m},q^{-x+N+1}\,;q)_k}{(q\,;q)_k}q^k.
  \label{PmoqH1}\!
\end{align}

\subsubsection{$\bm{q}$-Krawtchouk ($\bm{q}$K)}
\label{sec:qK1}

\begin{align}
  &\check{P}^{\text{monic}}_{N+1+m}(x;N,p)
  =\Lambda(x;N)\,q^{-(N+1)m}
  \check{P}^{\text{monic}}_m(x-N-1;-N-2,p'),
  \label{PmoqK}\\
  &\qquad p'=pq^{2(N+1)},\quad
  c_n(N,p)=\frac{(-pq^n\,;q)_n}{(q^{-N}\,;q)_n},\\
  &\check{P}^{\text{monic}}_m(x-N-1;-N-2,p')
  =\sum_{k=0}^m\frac{(q^{N+2+k}\,;q)_{m-k}}{(-pq^{2(N+1)+m+k}\,;q)_{m-k}}
  \frac{(q^{-m},q^{-x+N+1}\,;q)_k}{(q\,;q)_k}q^k.
  \label{PmoqK1}
\end{align}

\subsubsection{quantum $\bm{q}$-Krawtchouk (q$\bm{q}$K)}
\label{sec:qqK1}

\begin{align}
  &\check{P}^{\text{monic}}_{N+1+m}(x;N,p)
  =\Lambda(x;N)\,q^{-(N+1)m}
  \check{P}^{\text{monic}}_m(x-N-1;-N-2,p'),
  \label{PmoqqK}\\
  &\qquad p'=pq^{N+1},\quad
  c_n(N,p)=\frac{p^nq^{n^2}}{(q^{-N}\,;q)_n},\\
  &\check{P}^{\text{monic}}_m(x-N-1;-N-2,p')
  =\sum_{k=0}^m(q^{N+2+k}\,;q)_{m-k}
  \frac{(q^{-m},q^{-x+N+1}\,;q)_k}{(q\,;q)_k}\n[-4pt]
  &\phantom{\check{P}^{\text{monic}}_m(x-N-1;-N-2,p')=\sum_{k=0}^m}
  \times p^{k-m}q^{(N+m+2)k-m(m+N+1)}.
  \label{PmoqqK1}
\end{align}

\subsubsection{affine $\bm{q}$-Krawtchouk (a$\bm{q}$K)}
\label{sec:aqK1}

\begin{align}
  &\check{P}^{\text{monic}}_{N+1+m}(x;N,p)
  =\Lambda(x;N)\,q^{-(N+1)m}
  \check{P}^{\text{monic}}_m(x-N-1;-N-2,p').
  \label{PmoaqK}\\
  &\qquad p'=pq^{N+1},\quad c_n(N,p)=\frac{1}{(pq,q^{-N}\,;q)_n},\\
  &\check{P}^{\text{monic}}_m(x-N-1;-N-2,p')
  =\sum_{k=0}^m(pq^{N+2+k},q^{N+2+k}\,;q)_{m-k}
  \frac{(q^{-m},q^{-x+N+1}\,;q)_k}{(q\,;q)_k}q^k.
  \label{PmoaqK1}
\end{align}

\paragraph{(v) Family with $\bm{\eta(x)=(q^{-x}-1)(1-dq^x)}$}

\subsubsection{dual $\bm{q}$-Hahn (d$\bm{q}$H)}
\label{sec:dqH1}

\begin{align}
  &\check{P}^{\text{monic}}_{N+1+m}(x;N,\bm{\lambda})
  =\Lambda(x;N,\bm{\lambda})\,q^{-(N+1)m}
  \check{P}^{\text{monic}}_m(x-N-1;-N-2,\bm{\lambda}'),
  \label{PmodqH}\\
  &\qquad\qquad\bm{\lambda}'=(aq^{N+1},bq^{N+1}),\quad
  c_n(N,\bm{\lambda})=\frac{1}{(a,q^{-N}\,;q)_n},\\
  &\check{P}^{\text{monic}}_m(x-N-1;-N-2,\bm{\lambda}')
  =\sum_{k=0}^m(aq^{N+1+k},q^{N+2+k}\,;q)_{m-k}
  \frac{(q^{-m},abq^{x+N},q^{-x+N+1}\,;q)_k}{(q\,;q)_k}q^k.
  \label{PmodqH1}
\end{align}

\subsubsection{dual $\bm{q}$-Krawtchouk (d$\bm{q}$K)}
\label{sec:dqK1}

\begin{align}
  &\check{P}^{\text{monic}}_{N+1+m}(x;N,p)
  =\Lambda(x;N,p)\,q^{-(N+1)m}
  \check{P}^{\text{monic}}_m(x-N-1;-N-2,p'),
  \label{PmodqK}\\
  &\qquad p'=p^{2(N+1)},\quad c_n(N,p)=\frac{1}{(q^{-N}\,;q)_n},\\
  &\check{P}^{\text{monic}}_m(x-N-1;-N-2,p')
  =\sum_{k=0}^m(q^{N+2+k}\,;q)_{m-k}
  \frac{(q^{-m},q^{-x+N+1},-pq^{x+N+1}\,;q)_k}{(q\,;q)_k}q^k.
  \label{PmodqK1}
\end{align}

\subsection{Factorisation of multi-indexed and Krein-Adler systems}
\label{sec:set}

As shown in \S\,\ref{sec:fin}, the ``Diophantine'' property and factorisation
are a consequence of a general fact that the finite orthogonal polynomials in
the Askey scheme are the eigenvectors of certain real symmetric matrices
$\mathcal{H}$ \eqref{tdefcomp}. Certain generalisation of the finite
orthogonal polynomials in the Askey scheme are known for some time.
They are generated by multiple Darboux transformations \cite{crum} of these
polynomials by adopting certain seed solutions. When the polynomial themselves
are chosen as the seed solutions \cite{os22}, the obtained polynomials are
called Krein-Adler \cite{krein,adler} polynomials. When certain virtual state
vectors are used as the seed solutions, the new polynomials are called,
for example, multi-indexed ($q$-)Racah polynomials, etc \cite{os26}.
By construction, {\em i.e.} by appropriate choices of the seed solutions and
the parameter ranges, these new types of orthogonal polynomials are the
eigenvectors of certain real symmetric matrices. Therefore, when
$\check{P}_n(x)$ is changed to $\check{P}^{\text{monic}}_{N+1+m}(x)$ in
\eqref{phiMform}, the new types of orthogonal polynomials have the
``Diophantine'' property and factorisation.

As will be discussed in the subsequent two sections, a new type of
multi-indexed polynomials can be generated by using
$\check{P}^{\text{monic}}_{N+1+m}(x)$ as seed solutions for the multiple
Darboux transformations. If the resulting orthogonality measures are positive
definite, which is the case for those introduced in \S\,\ref{sec:shape},
such systems have also the ``Diophantine'' property and factorisation.

\section{New Multi-indexed Polynomials, General Setting}
\label{sec:newmulti}

One possible application of the new polynomials satisfying the ``Diophantine''
and factorisation properties \eqref{Dfac} is the generation of new
multi-indexed orthogonal polynomials. It is well known that certain infinite
norm solutions of the Schr\"odinger equations with the radial oscillator and
P\"oschl-Teller potentials are used as seed solutions for multiple Darboux
transformations to generate multi-indexed Laguerre and Jacobi polynomials
\cite{gomez,quesne,os16,os25}.
Likewise the non-eigenvector solutions of the tri-diagonal real symmetric
matrix $\mathcal{H}$ eigenvalue problem \eqref{comp} are used as seed
solutions to generate the multi-indexed versions of the ($q$)-Racah
polynomials \cite{os26} etc through the difference equation analogues of the
multiple Darboux transformations.
Now we have plenty of zero norm solutions
\begin{align}
  \tilde{\phi}_m(x)&=\phi^{\text{monic}}_{N+1+m}(x)
  =\phi_0(x)\check{P}^{\text{monic}}_{N+1+m}(x;N,\bm{\lambda})
  \quad(m\in\mathbb{Z}_{\ge0})
  \label{zero1}\\
  &=\phi_0(x)\Lambda(x;N,\bm{\lambda})\check{Q}_m(x;N,\bm{\lambda}),\quad
  \bigl(\tilde{\phi}_m(x)\bigr)_{x\in\cX}\in\mathbb{R}^{N+1},
  \label{zero2}
\end{align}
corresponding to \eqref{Pmonic4}, \eqref{PmoKr}, \eqref{PmoHa}, \eqref{PmoRa},
\eqref{PmodH}, \eqref{PmodqqK}, \eqref{PmoqH}, \eqref{PmoqK}, \eqref{PmoqqK},
\eqref{PmoaqK}, \eqref{PmodqH} and \eqref{PmodqK},
which could be used as seed solutions to generate the multi-indexed versions
of the twelve polynomials listed in \S\,\ref{sec:Kra}--\S\,\ref{sec:dqK}.

Let us recapitulate the basic formulas of the multiple Darboux transformation
of the original eigenvectors $\{\phi_n(x)\}$ \eqref{phindef} by using the
zero-norm vectors $\{\tilde{\phi}_m(x)\}$ \eqref{zero1} as seed solutions.
These formulas apply for each of the finite polynomials listed in
\S\,\ref{sec:data}. We employ the formulas reported in
\S\,2.3 of \cite{os24}.

We choose $M$ distinct zero-norm seed solutions
\begin{equation}
  \{\tilde{\phi}_{m_1}(x),\tilde{\phi}_{m_2}(x),\ldots,
  \tilde{\phi}_{m_M}(x)\},\quad 
  0\le m_1< m_2 <\cdots < m_M,
  \label{seedsols}
\end{equation}
which correspond to the index set
\begin{equation}
  \mathcal{D}=\{m_1,\ldots,m_M\}.
\end{equation}
The deformed set of eigenvectors $\{\bar{\phi}_{\mathcal{D},n}\}$ ($n\in\cX$)
satisfy the difference equation
\begin{align}
  &\mathcal{H}_{\mathcal{D}}\bar{\phi}_{\mathcal{D},n}(x)
  =\mathcal{E}(n)\bar{\phi}_{\mathcal{D},n}(x),\\
  &\mathcal{H}_{\mathcal{D}}\eqdef
  \mathcal{A}^\dagger_{\mathcal{D}}\mathcal{A}_{\mathcal{D}},\quad
  \mathcal{A}_{\mathcal{D}}\eqdef\sqrt{\bar{B}(x)}-e^\partial\sqrt{\bar{D}(x)},
  \quad \mathcal{A}_{\mathcal{D}}^\dagger\eqdef
  \sqrt{\bar{B}(x)}-\sqrt{\bar{D}(x)}e^{-\partial},
  \label{dbfQMAAd}\\
  &(\bar{\phi}_{\mathcal{D},n},\bar{\phi}_{\mathcal{D},\ell})
  =\prod_{j=1}^{M}\bigl(\mathcal{E}(n)-\mathcal{E}(N+1+m_j)\bigr)
  \cdot\frac{\delta_{n\,\ell}}{d_n^2}.
  \label{modnorm}
\end{align}
In these formulas $e^\partial$ ($\partial\eqdef\frac{d}{dx}$) is a finite
shift operator acting on smooth functions,
\begin{equation}
  e^{\partial}f(x)=f(x+1),\quad(e^{\partial})^j=e^{j\partial},\quad
  e^{j\partial}f(x)=f(x+j)\quad(j\in\mathbb{Z}).
  \label{finshift}
\end{equation}
As an operator it passes a smooth function $g(x)$ as
\begin{equation*}
  e^{j\partial}g(x)=g(x+j)e^{j\partial}.
\end{equation*}
Here $\bar{B}(x)$, $\bar{D}(x)$ and $\bar{\phi}_{\mathcal{D},n}(x)$ are
expressed in terms of the Casoratians involving the seed solutions
$\{\tilde{\phi}_m\}$, eigenvector $\phi_n$ and the groundstate eigenvector
$\phi_0$,
\begin{align}
  &\text{W}_{\text{C}}[f_1,\ldots,f_n](x)
  \eqdef\det\Bigl(f_k(x+j-1)\Bigr)_{1\leq j,k\leq n},
  \label{Casdef}\\
  &\bar{B}(x)=\sqrt{B(x+M)D(x+M+1)}\,
  \frac{\text{W}_{\text{C}}[\tilde{\phi}_{m_1},\ldots,\tilde{\phi}_{m_{M}}](x)}
  {\text{W}_{\text{C}}[\tilde{\phi}_{m_1},\ldots,\tilde{\phi}_{m_{M}}](x+1)}\n
  &\phantom{\bar{B}(x)=}\times
  \frac{\text{W}_{\text{C}}[\tilde{\phi}_{m_1},\ldots,\tilde{\phi}_{m_{M}},
  {\phi}_{0}](x+1)}
  {\text{W}_{\text{C}}[\tilde{\phi}_{m_1},\ldots,\tilde{\phi}_{m_{M}},
  \phi_{0}](x)},
  \label{Bbar}\\
  &\bar{D}(x)=\sqrt{B(x-1)D(x)}\,
  \frac{\text{W}_{\text{C}}[\tilde{\phi}_{m_1},\ldots,\tilde{\phi}_{m_{M}}](x+1)}
  {\text{W}_{\text{C}}[\tilde{\phi}_{m_1},\ldots,\tilde{\phi}_{m_{M}}](x)}\n
  &\phantom{\bar{D}(x)=}\times
  \frac{\text{W}_{\text{C}}[\tilde{\phi}_{m_1},\ldots,\tilde{\phi}_{m_{M}},
  \phi_{0}](x-1)}
  {\text{W}_{\text{C}}[\tilde{\phi}_{m_1},\ldots,\tilde{\phi}_{m_{M}},
  \phi_{0}](x)},
  \label{Dbar}\\
  &F(x)\eqdef\frac{\sqrt{\prod_{k=1}^{M}B(x+k-1)D(x+k)}}
  {\text{W}_{\text{C}}[\tilde{\phi}_{m_1},\ldots,\tilde{\phi}_{m_{M}}](x)\,
  \text{W}_{\text{C}}[\tilde{\phi}_{m_1},\ldots,\tilde{\phi}_{m_{M}}](x+1)},
  \n[2pt]
  &\bar{\phi}_{\mathcal{D},n}(x)=
  (-1)^{M}\sqrt{F(x)}\,
  \text{W}_{\text{C}}[\tilde{\phi}_{m_1},\ldots,\tilde{\phi}_{m_{M}},\phi_n](x).
  \label{phibarnD}
\end{align}
By using the Casoratian identity
\begin{align}
  &\text{W}_{\text{C}}[gf_1,gf_2,\ldots,gf_n](x)
  =\prod_{k=0}^{n-1}g(x+k)\cdot\text{W}_{\text{C}}[f_1,f_2,\ldots,f_n](x),
  \label{Wformula1}
\end{align}
some formulas are simplified,
\begin{align*} 
  \frac{\text{W}_{\text{C}}[\tilde{\phi}_{m_1},\ldots,\tilde{\phi}_{m_{M}}](x)}
  {\text{W}_{\text{C}}[\tilde{\phi}_{m_1},\ldots,\tilde{\phi}_{m_{M}}](x+1)}
  &=\frac{\phi_0(x)\Lambda(x)}{\phi_0(x+M)\Lambda(x+M)}
  \frac{\text{W}_{\text{C}}[\check{Q}_{m_1},\ldots,\check{Q}_{m_{M}}](x)}
  {\text{W}_{\text{C}}[\check{Q}_{m_1},\ldots,\check{Q}_{m_{M}}](x+1)},\\
  \frac{\text{W}_{\text{C}}[\tilde{\phi}_{m_1},\ldots,\tilde{\phi}_{m_{M}}](x+1)}
  {\text{W}_{\text{C}}[\tilde{\phi}_{m_1},\ldots,\tilde{\phi}_{m_{M}}](x)}
  &=\frac{\phi_0(x+M)\Lambda(x+M)}{\phi_0(x)\Lambda(x)}
  \frac{\text{W}_{\text{C}}[\check{Q}_{m_1},\ldots,\check{Q}_{m_{M}}](x+1)}
  {\text{W}_{\text{C}}[\check{Q}_{m_1},\ldots,\check{Q}_{m_{M}}](x)},\\
  \frac{\text{W}_{\text{C}}[\tilde{\phi}_{m_1},\ldots,\tilde{\phi}_{m_{M}},
  {\phi}_{0}](x+1)}
  {\text{W}_{\text{C}}[\tilde{\phi}_{m_1},\ldots,\tilde{\phi}_{m_{M}},
  \phi_{0}](x)}
  &=\frac{\phi_0(x+M+1)\Lambda(x+M+1)}{\phi_0(x)\Lambda(x)}\n
  &\quad\times
  \frac{\text{W}_{\text{C}}[\check{Q}_{m_1},\ldots,\check{Q}_{m_{M}},
  \Lambda^{-1}](x+1)}
  {\text{W}_{\text{C}}[\check{Q}_{m_1},\ldots,\check{Q}_{m_{M}},
  \Lambda^{-1}](x)},\\
  \frac{\text{W}_{\text{C}}[\tilde{\phi}_{m_1},\ldots,\tilde{\phi}_{m_{M}},
  {\phi}_{0}](x-1)}
  {\text{W}_{\text{C}}[\tilde{\phi}_{m_1},\ldots,\tilde{\phi}_{m_{M}},
  \phi_{0}](x)}
  &=\frac{\phi_0(x-1)\Lambda(x-1)}{\phi_0(x+M)\Lambda(x+M)}
  \frac{\text{W}_{\text{C}}[\check{Q}_{m_1},\ldots,\check{Q}_{m_{M}},
  \Lambda^{-1}](x-1)}
  {\text{W}_{\text{C}}[\check{Q}_{m_1},\ldots,\check{Q}_{m_{M}},
  \Lambda^{-1}](x)},\\
  {\text{W}_{\text{C}}[\tilde{\phi}_{m_1},\ldots,\tilde{\phi}_{m_{M}},
  \phi_{n}](x)}
  &=\prod_{k=0}^M\phi_0(x+k)\Lambda(x+k)\cdot
  {\text{W}_{\text{C}}[\check{Q}_{m_1},\ldots,\check{Q}_{m_{M}},
  \Lambda^{-1}\check{P}_n](x)}.
\end{align*}
These lead to
\begin{align}
  \bar{B}(x)&=\sqrt{B(x+M)D(x+M+1)}\,
  \frac{\phi_0(x+M+1)\Lambda(x+M+1)}{\phi_0(x+M)\Lambda(x+M)}\n
  &\quad\times
  \frac{\text{W}_{\text{C}}[\check{Q}_{m_1},\ldots,\check{Q}_{m_{M}}](x)}
  {\text{W}_{\text{C}}[\check{Q}_{m_1},\ldots,\check{Q}_{m_{M}}](x+1)}\,
  \frac{\text{W}_{\text{C}}[\check{Q}_{m_1},\ldots,\check{Q}_{m_{M}},
  \Lambda^{-1}](x+1)}
  {\text{W}_{\text{C}}[\check{Q}_{m_1},\ldots,\check{Q}_{m_{M}},
  \Lambda^{-1}](x)}\n
  &=B(x+M)\frac{\Lambda(x+M+1)}{\Lambda(x+M)}\n
  &\quad\times
  \frac{\text{W}_{\text{C}}[\check{Q}_{m_1},\ldots,\check{Q}_{m_{M}}](x)}
  {\text{W}_{\text{C}}[\check{Q}_{m_1},\ldots,\check{Q}_{m_{M}}](x+1)}\,
  \frac{\text{W}_{\text{C}}[\check{Q}_{m_1},\ldots,\check{Q}_{m_{M}},
  \Lambda^{-1}](x+1)}
  {\text{W}_{\text{C}}[\check{Q}_{m_1},\ldots,\check{Q}_{m_{M}},
  \Lambda^{-1}](x)},
  \label{BMform}\\
  \bar{D}(x)&=\sqrt{B(x-1)D(x)}\,
  \frac{\phi_0(x-1)\Lambda(x-1)}{\phi_0(x)\Lambda(x)}\n
  &\quad\times
  \frac{\text{W}_{\text{C}}[\check{Q}_{m_1},\ldots,\check{Q}_{m_{M}}](x+1)}
  {\text{W}_{\text{C}}[\check{Q}_{m_1},\ldots,\check{Q}_{m_{M}}](x)}
  \frac{\text{W}_{\text{C}}[\check{Q}_{m_1},\ldots,\check{Q}_{m_{M}},
  \Lambda^{-1}](x-1)}
  {\text{W}_{\text{C}}[\check{Q}_{m_1},\ldots,\check{Q}_{m_{M}},
  \Lambda^{-1}](x)}\n
  &=D(x)\frac{\Lambda(x-1)}{\Lambda(x)}\n
  &\quad\times
  \frac{\text{W}_{\text{C}}[\check{Q}_{m_1},\ldots,\check{Q}_{m_{M}}](x+1)}
  {\text{W}_{\text{C}}[\check{Q}_{m_1},\ldots,\check{Q}_{m_{M}}](x)}
  \frac{\text{W}_{\text{C}}[\check{Q}_{m_1},\ldots,\check{Q}_{m_{M}},
  \Lambda^{-1}](x-1)}
  {\text{W}_{\text{C}}[\check{Q}_{m_1},\ldots,\check{Q}_{m_{M}},
  \Lambda^{-1}](x)}.
  \label{DMform}
\end{align}
Furthermore, we have
\begin{align*}
  &\quad
  {\text{W}_{\text{C}}[\tilde{\phi}_{m_1},\ldots,\tilde{\phi}_{m_{M}}](x)\,
  \text{W}_{\text{C}}[\tilde{\phi}_{m_1},\ldots,\tilde{\phi}_{m_{M}}](x+1)}\n
  &=\bigl(\phi_0(x)\phi_0(x+M)\Lambda(x)\Lambda(x+M)\bigr)^{-1}
  \prod_{k=0}^M\phi_0(x+k)^2\Lambda(x+k)^2\n
  &\quad\times{\text{W}_{\text{C}}[\check{Q}_{m_1},\ldots,\check{Q}_{m_{M}}](x)}
  {\text{W}_{\text{C}}[\check{Q}_{m_1},\ldots,\check{Q}_{m_{M}}](x+1)},\\
  &\sqrt{\prod_{k=1}^{M}B(x+k-1)D(x+k)}\ \phi_0(x)\phi_0(x+M)
  =\phi_0(x)^2\prod_{k=1}^{M}B(x+k-1).
\end{align*}
In terms of these we arrive at
\begin{align}
  F(x)&=\prod_{k=1}^{M}B(x+k-1)\cdot\frac{\Lambda(x+M)}{\Lambda(x)}
  \frac1{\prod_{k=1}^M\phi_0(x+k)^2\Lambda(x+k)^2}\n
  &\quad\times
  \frac{1}{\text{W}_{\text{C}}[\check{Q}_{m_1},\ldots,\check{Q}_{m_{M}}](x)
  \text{W}_{\text{C}}[\check{Q}_{m_1},\ldots,\check{Q}_{m_{M}}](x+1)},
  \label{FMform}\\
  \bar{\phi}_{\mathcal{D},n}(x)&=
  (-1)^{M}\sqrt{\prod_{k=1}^{M}B(x+k-1)}\ \phi_0(x)
  \sqrt{\Lambda(x)\Lambda(x+M)}\n
  &\quad\times
  \frac{\text{W}_{\text{C}}[\check{Q}_{m_1},\ldots,\check{Q}_{m_{M}},
  \Lambda^{-1}\check{P}_n](x)}
  {\sqrt{\text{W}_{\text{C}}[\check{Q}_{m_1},\ldots,\check{Q}_{m_{M}}](x)
  \text{W}_{\text{C}}[\check{Q}_{m_1},\ldots,\check{Q}_{m_{M}}](x+1)}}.
  \label{phiMform}
\end{align}
 
For each type of the polynomials, the multi-indexed polynomials are identified
by removing various kinematical factors from $\bar{\phi}_{\mathcal{D},n}$.
As shown by \eqref{modnorm}, the orthogonality is built in. After the
identification and the orthogonality, securing the positivity of the
orthogonality measures is essential. It must be verified for each polynomial,
for each choice of the set $\mathcal{D}$ of the zero-norm seed solutions and
for the appropriate ranges of the involved parameters.
The situation is much more complicated than the multi-indexed Laguerre and
Jacobi polynomials \cite{os25} and ($q$)-Racah polynomials \cite{os26}, etc.
In these established cases, the seed solutions can be chosen to have a
definite sign, and that is closely related to the positive definite
orthogonality measures.
We hope these detailed tasks for each type of polynomials will be carried out
in future.
A few examples of single-indexed exceptional Krawtchouk polynomials are
reported in \cite{miki} and their ``Diophantine'' property and factorisation
are discussed.

\section{Shape-invariant Cases}
\label{sec:shape}

Here we report on multi-indexed polynomials corresponding to a very special
choice of the seed solutions, $M$ contiguous lowest degree zero-norm solutions,
\begin{equation}
  \mathcal{D}=\{0,1,\ldots,M-1\}.
  \label{0M}
\end{equation}
For the Schr\"odinger equations, for which the Laguerre and Jacobi polynomials
are the eigen polynomials, the multiple Darboux transformations by using $M$
contiguous lowest degree eigenfunctions produce remarkable effects
\cite{crum,krein,adler},
\begin{alignat*}{3}
  \text{Laguerre}:&\ \ &L_n^{(\alpha)}\bigl(\eta(x)\bigr)
  &\,\to\,(\text{const})\times L_{n-M}^{(\alpha+M)}\bigl(\eta(x)\bigr)
  &\ \ (n\ge M),\\
  \text{Jacobi}:&\ \ &P_n^{(\alpha,\beta)}\bigl(\eta(x)\bigr)
  &\,\to\,(\text{const})\times P_{n-M}^{(\alpha+M,\beta+M)}\bigl(\eta(x)\bigr)
  &\ \ (n\ge M).
\end{alignat*}
Under the transformations, the polynomials keep their identity with shifted
parameters and the positivity of the orthogonality measure unchanged.
This phenomenon is called shape-invariance \cite{genden}.
Shape-invariance also holds for all the polynomials listed in \S\,\ref{sec:data}
and the others in the Askey scheme. It is called the forward shift relation
\cite{kls}. For example
\begin{alignat*}{3}
  \text{Racah}:&\ \ &\check{P}_n(x;N,b,c,d)
  &\,\to\,(\text{const})\times\check{P}_{n-M}(x;N-M,b+M,c+M,d+M)&\ \ (n\ge M),\\
  \text{$q$-Racah}:&\ \ &\check{P}_n(x;N,b,c,d)
  &\,\to\,(\text{const})\times\check{P}_{n-M}(x;N-M,bq^M,cq^M,dq^M)&\ \ (n\ge M).
\end{alignat*}
For more details, see \cite{os12}.

We will demonstrate similar effects for the present case of using $M$
contiguous lowest zero-norm states \eqref{0M} for the finite orthogonal
polynomials.
The explicit formulas of the transformation of $\bar{B}$ \eqref{BMform},
$\bar{D}$ \eqref{DMform} and $ \bar{\phi}_{\mathcal{D},n}$ \eqref{phiMform},
depend on the type of the sinusoidal coordinates $\eta(x)$.
For three families (i), (iii) and (iv), in which the sinusoidal coordinate
contains no parameter other than $q$, (i) $\eta(x)=x$, (iii) $\eta(x)=1-q^x$
and (iv) $\eta(x)=q^{-x}-1$.
For the other two families (ii) and (v), in which $\eta(x)$ contains $d$,
(ii) $\eta(x)=x(x+d)$, (v) $\eta(x)=(q^{-x}-1)(1-dq^x)$, the parameter $d$
shifts in the formulas. They are summarised in the following
\begin{theo}
\label{theo:sim}
The results of the multiple Darboux transformation are
\begin{alignat}{2}
  \text{Family {\rm (i),\,(iv)}}:&&\quad B(x;N,\bm{\lambda})&
  \,\to\,B(x+M;N+M,\bm{\lambda}),
  \label{f14B}\\
  &&D(x;N,\bm{\lambda})&\,\to\,D(x+M;N+M,\bm{\lambda}),
  \label{f14D}\\
  &&\check{P}_n(x;N,\bm{\lambda})&\,\to\,
  (\text{\rm const})\times\check{P}_n(x+M;N+M,\bm{\lambda}),
  \label{f14P}\\
  \text{Family {\rm (ii)}}:&&\quad B(x;N,d,\bar{\bm{\lambda}})
  &\,\to\,B(x+M;N+M,d-M,\bar{\bm{\lambda}}),
  \label{f2B}\\
  &&D(x;N,d,\bar{\bm{\lambda}})&\,\to\,D(x+M;N+M,d-M,\bar{\bm{\lambda}}),
  \label{f2D}\\
  &&\check{P}_n(x;N,d,\bar{\bm{\lambda}})&\,\to\,
  (\text{\rm const})\times\check{P}_n(x+M;N+M,d-M,\bar{\bm{\lambda}}),
  \label{f2P}\\
  \text{Family {\rm (iii)}}:&&\quad B(x;N,p)
  &\,\to\,B(x+M;N+M,pq^{-M}),
  \label{f3B}\\
  &&D(x;N,p)&\,\to\,D(x+M;N+M,pq^{-M}),
  \label{f3D}\\
  &&\check{P}_n(x;N,p)&\,\to\,
  (\text{\rm const})\times\check{P}_n(x+M;N+M,pq^{-M}),
  \label{f3P}\\
  \text{Family {\rm (v)}}:&&\quad B(x;N,d,\bar{\bm{\lambda}})
  &\,\to\,B(x+M;N+M,dq^{-M},\bar{\bm{\lambda}}),
  \label{f5B}\\
  &&D(x;N,d,\bar{\bm{\lambda}})&\,\to\,D(x+M;N+M,dq^{-M},\bar{\bm{\lambda}}),
  \label{f5D}\\
  &&\check{P}_n(x;N,d,\bar{\bm{\lambda}})&\,\to\,
  (\text{\rm const})\times\check{P}_n(x+M;N+M,dq^{-M},\bar{\bm{\lambda}}),
  \label{f5P}
\end{alignat}
in which $n\in\cX$ and $\bar{\bm{\lambda}}$ stands for the parameters other
than $N$ and $d$.
The positivity of the orthogonality measures is unchanged for the formulas
\eqref{f14B}--\eqref{f14P}. The positivity also holds for the formulas
\eqref{f2B}--\eqref{f5P} for an appropriate range of the parameter $p$ or $d$.
\end{theo}

In the rest of this section we present the outline of the derivation of the
formulas \eqref{f14B}--\eqref{f5P}.

In order to obtain the explicit forms for $\bar{B}(x)$ \eqref{BMform},
$\bar{D}(x)$ \eqref{DMform} and $\bar{\phi}_{\mathcal{D},n}$ \eqref{phiMform}
we need to evaluate only four formulas
\begin{align*}
  \text{W}_{\text{C}}[1,\check{Q}_{1},\ldots,\check{Q}_{M-1}](x)&,
  &\hspace{-5mm}&
  \Lambda(x){\text{W}_{\text{C}}[1,\check{Q}_{1},\ldots,\check{Q}_{M-1},
  \Lambda^{-1}](x)},\\
  \Lambda(x+M){\text{W}_{\text{C}}[1,\check{Q}_{1},\ldots,\check{Q}_{M-1},
  \Lambda^{-1}](x)}&,
  &\hspace{-5mm}&
  \Lambda(x)\text{W}_{\text{C}}[1,\check{Q}_{1},\ldots,\check{Q}_{M-1},
  \Lambda^{-1}\check{P}_n](x).
\end{align*}
The other components are obtained by shifting $x\to x\pm1$.
As shown in {\bf Theorem\,\ref{theo1}}, $\check{Q}_m(x;N,\bm{\lambda})$ is a
monic degree $m$ polynomial in $\eta(x;\bm{\lambda})$. Therefore these
formulas are simplified by sweeping the Casoratian determinants successively,
\begin{align}
  \text{W}_{\text{C}}[1,\check{Q}_{1},\ldots,\check{Q}_{M-1}](x)&=
  \text{W}_{\text{C}}[1,\eta,\eta^2,\ldots,\eta^{M-1}](x),
  \label{vand1}\\
  \Lambda(x){\text{W}_{\text{C}}[1,\check{Q}_{1},\ldots,\check{Q}_{M-1},
  \Lambda^{-1}](x)}
  &=\Lambda(x)\text{W}_{\text{C}}[1,\eta,\eta^2,\ldots,\eta^{M-1},
  \Lambda^{-1}](x),
  \label{vand2}\\
  \Lambda(x+M){\text{W}_{\text{C}}[1,\check{Q}_{1},\ldots,\check{Q}_{M-1},
  \Lambda^{-1}](x)}
  &=\Lambda(x+M)\text{W}_{\text{C}}[1,\eta,\eta^2,\ldots,\eta^{M-1},
  \Lambda^{-1}](x),
  \label{vand3}\\
  \Lambda(x)\text{W}_{\text{C}}[1,\check{Q}_{1},\ldots,\check{Q}_{M-1},
  \Lambda^{-1}\check{P}_n](x)
  &=\Lambda(x)\text{W}_{\text{C}}[1,\eta,\eta^2,\ldots,\eta^{M-1},
  \Lambda^{-1}\check{P}_n](x).
  \label{vand4}
\end{align}
It should be stressed that the explicit expressions of $\check{Q}_m(x)$ listed
in \S\,\ref{sec:Kra}--\S\,\ref{sec:dqK} are not necessary. Now that the
r.h.s.\ of \eqref{vand1} is explicitly known as it is a Vandermonde determinant.
Likewise \eqref{vand2} and \eqref{vand3} are easily evaluated as
$\Lambda(x)/\Lambda(x+j-1)$ and $\Lambda(x+M)/\Lambda(x+j-1)$ have simple
expressions. It is important to stress that these quantities depend on $N$
and $\eta(x;\bm{\lambda})$ only, except for the one \eqref{vand4}
containing $\check{P}_n$.

Now we list the explicit forms of the expressions of
\eqref{vand1}--\eqref{vand4} for each family (i)--(v).
By using the following constants,
\begin{align}
  c(M)&\eqdef\prod_{k=1}^{M-1}k!,\quad
  c(N,M)\eqdef(-1)^M\,(N+1)_M\,c(M),\\
  c_q(M)&\eqdef q^{-\frac16M(M-1)(2M-1)}\prod_{k=1}^{M-1}(q\,;q)_k,\\
  c_q(N,M)&\eqdef (-1)^Mq^{-\frac12M(M-1)}(q^{N+1}\,;q)_M\,c_q(M),\\
  \binom{M}{j}&=\frac{M!}{j!\,(M-j)!},\quad
  \genfrac{[}{]}{0pt}{}{\,M\,}{j}\eqdef
  \frac{(q\,;q)_M}{(q\,;q)_j(q\,;q)_{M-j}},
\end{align}
they are given as follows.\\
Family (i) :
\begin{align}
  &\text{W}_{\text{C}}[1,x,x^2,\ldots,x^{M-1}](x)=c(M),\\
  &\Lambda(x)\text{W}_{\text{C}}[1,x,x^2,\ldots,x^{M-1},\Lambda^{-1}](x)
  =\frac{c(N,M)}{(x+1)_M},\\
  &\Lambda(x+M)\text{W}_{\text{C}}[1,x,x^2,\ldots,x^{M-1},\Lambda^{-1}](x)
  =\frac{c(N,M)}{(x-N)_M},\\
  &\Lambda(x)
  \text{W}_{\text{C}}[1,x,x^2,\ldots,x^{M-1},\Lambda^{-1}\check{P}_n](x)\n
  &=\frac{(-1)^Mc(M)}{(x+1)_M}
  \sum_{j=0}^M(-1)^j\binom{M}{j}
  (x+1+j)_{M-j}(x-N)_j\check{P}_n(x+j),
  \label{1eqform}
\end{align}
Family (ii) :
\begin{align}
  &\text{W}_{\text{C}}[1,\eta,\eta^2,\ldots,\eta^{M-1}](x)=
  c(M)\prod_{k=1}^{M-1}(2x+k+d)_k,\\
  &\Lambda(x)\text{W}_{\text{C}}[1,\eta,\eta^2,\ldots,\eta^{M-1},\Lambda^{-1}](x)
  =\frac{c(N,M)\prod_{k=1}^M(2x+k+d)_k}{(x+1,x+N+1+d)_M},\\
  &\Lambda(x+M)
  \text{W}_{\text{C}}[1,\eta,\eta^2,\ldots,\eta^{M-1},\Lambda^{-1}](x)
  =\frac{c(N,M)\prod_{k=1}^M(2x+k+d)_k}{(x-N,x+d)_M},\\
  &\Lambda(x)\text{W}_{\text{C}}[1,\eta,\eta^2,\ldots,\eta^{M-1},
  \Lambda^{-1}\check{P}_n](x)\n
  &=\frac{(-1)^Mc(M)\prod_{k=1}^{M-2}(2x+2+k+d)_k}{(x+1,x+N+1+d)_M}
  \sum_{j=0}^M(-1)^j\binom{M}{j}\mathcal{T}(x,M,j,d)\n
  &\qquad\times
  (x+1+j,x+1+j+N+d)_{M-j}(x-N,x+d)_j\check{P}_n(x+j),
  \label{2eqform}\\
  &\mathcal{T}(x,M,j,d)\eqdef(2x+M+1+j+d)_{M-j-1}(2x+1+d)_{j-1}\n
  &\phantom{\mathcal{T}(x,M,j,d)\eqdef}\times\left\{
  \begin{array}{cl}
  1&:j=0,M\\
  (2x+2j+d)&:\text{otherwise}
  \end{array}\right.,
  \label{Fdef}
\end{align}
Family (iii) :
\begin{align}
  &\text{W}_{\text{C}}[1,\eta,\eta^2,\ldots,\eta^{M-1}](x)
  =c_q(M)\,q^{\frac12M(M-1)x}q^{\frac12M(M-1)^2},\\
  &\Lambda(x)\text{W}_{\text{C}}[1,\eta,\eta^2,\ldots,\eta^{M-1},\Lambda^{-1}](x)
  =\frac{c_q(N,M)\,q^{\frac12M(M+1)x}q^{\frac12M(M^2-2N-1)}}{(q^{x+1}\,;q)_M},\\
  &\Lambda(x+M)
  \text{W}_{\text{C}}[1,\eta,\eta^2,\ldots,\eta^{M-1},\Lambda^{-1}](x)
  =\frac{c_q(N,M)\,q^{\frac12M(M+1)x}q^{\frac12M(M^2-2N-1)}}{(q^{x-N}\,;q)_M},\\
  &\Lambda(x)\text{W}_{\text{C}}[1,\eta,\eta^2,\ldots,\eta^{M-1},
  \Lambda^{-1}\check{P}_n](x)\n
  &=\frac{(-1)^Mc_q(M)q^{\frac12M(M-1)x}q^{\frac12M^2(M-1)-MN}}{(q^{x+1}\,q)_M}
  \sum_{j=0}^M(-1)^j\genfrac{[}{]}{0pt}{}{\,M\,}{j}q^{\frac12j(j+1)+M(N-j)}\n
  &\hspace{80mm}\times
  (q^{x+1+j}\,;q)_{M-j}(q^{x-N}\,;q)_j\check{P}_n(x+j),
  \label{3eqform}
\end{align}
Family (iv) :
\begin{align}
  &\text{W}_{\text{C}}[1,\eta,\eta^2,\ldots,\eta^{M-1}](x)=
  c_q(M)\,q^{-\frac12M(M-1)x},\\
  &\Lambda(x)\text{W}_{\text{C}}[1,\eta,\eta^2,\ldots,\eta^{M-1},\Lambda^{-1}](x)
  =\frac{c_q(N,M)\,q^{-\frac12M(M-1)x}}{(q^{x+1}\,;q)_M},\\
  &\Lambda(x+M)
  \text{W}_{\text{C}}[1,\eta,\eta^2,\ldots,\eta^{M-1},\Lambda^{-1}](x)
  =\frac{c_q(N,M)\,q^{-\frac12M(M-1)x}}{q^{M(N+1)}(q^{x-N}\,;q)_M},\\
  &\Lambda(x)\text{W}_{\text{C}}[1,\eta,\eta^2,\ldots,\eta^{M-1},
  \Lambda^{-1}\check{P}_n](x)\n
  &=\frac{(-1)^Mc_q(M)\,q^{-\frac12M(M-1)x}q^{-\frac12M(M-1)}}{(q^{x+1}\,;q)_M}
  \sum_{j=0}^M(-1)^j\genfrac{[}{]}{0pt}{}{\,M\,}{j}
  q^{\frac12j(j+1)+Nj}\n
  &\hspace{75mm}\times
  (q^{x+1+j}\,;q)_{M-j}(q^{x-N}\,;q)_j\check{P}_n(x+j),
  \label{4eqform}
\end{align}
Family (v) :
\begin{align}
  &\text{W}_{\text{C}}[1,\eta,\eta^2,\ldots,\eta^{M-1}](x)=
  c_q(M)\,q^{-\frac12M(M-1)x}\prod_{k=1}^{M-1}(dq^{2x+k}\,;q)_k,\\
  &\Lambda(x)\text{W}_{\text{C}}[1,\eta,\eta^2,\ldots,\eta^{M-1},\Lambda^{-1}](x)
  =\frac{c_q(N,M)\,q^{-\frac12M(M-1)x}
  \prod_{k=1}^M(dq^{2x+k}\,;q)_k}{(q^{x+1},dq^{x+1+N}\,;q)_M},\\
  &\Lambda(x+M)
  \text{W}_{\text{C}}[1,\eta,\eta^2,\ldots,\eta^{M-1},\Lambda^{-1}](x)
  =\frac{c_q(N,M)\,q^{-\frac12M(M-1)x}
  \prod_{k=1}^M(dq^{2x+k}\,;q)_k}{q^{M(N+1)}(q^{x-N},dq^x\,;q)_M},\\[4pt]
  &\Lambda(x)\text{W}_{\text{C}}[1,\eta,\eta^2,\ldots,\eta^{M-1},
  \Lambda^{-1}\check{P}_n](x)\n
  &=\frac{(-1)^Mc_q(M)\,q^{-\frac12M(M-1)x}q^{-\frac12M(M-1)}
  \prod_{k=1}^{M-2}(dq^{2x+2+k}\,;q)_k}
  {(q^{x+1},dq^{x+N+1}\,;q)_M}
  \sum_{j=0}^M(-1)^j\genfrac{[}{]}{0pt}{}{\,M\,}{j}
  q^{\frac12j(j+1)+Nj}\n
  &\qquad\times\mathcal{T}_q(x,M,j,d)(q^{x+1+j},dq^{x+1+j+N}\,;q)_{M-j}
  (q^{x-N},dq^x\,;q)_j\check{P}_n(x+j),
  \label{5eqform}\\
  &\mathcal{T}_q(x,M,j,d)\eqdef(dq^{2x+M+1+j}\,;q)_{M-j-1}(dq^{2x+1}\,;q)_{j-1}\,
  \times\left\{
  \begin{array}{cl}
  1&:j=0,M\\
  (1-dq^{2x+2j})&:\text{otherwise}
  \end{array}\right.\!.
  \label{qFdef}
\end{align}
Based on these results, the transformed functions $\bar{B}$ \eqref{BMform} and 
$\bar{D}$ \eqref{DMform} are calculated for each family:
\begin{align}
  \text{Family (i)}:\ \ \bar{B}(x)
  &=B(x+M)\,\frac{(x-N)_M}{(x+1-N)_M}=B(x+M)\,\frac{N-x}{N-x-M},\\
  \bar{D}(x)&=D(x)\,\frac{(x+1)_M}{(x)_M}=D(x)\,\frac{x+M}{x},\\
  \text{Family (ii)}:\ \ \bar{B}(x)
  &=B(x+M)\,\frac{(N-x)(x+d)(2x+2M+d)_2}{(N-x-M)(x+d+M)(2x+M+d)_2},\\
  \bar{D}(x)&=D(x)\,\frac{(x+M)(x+N+M+d)(2x-1+d)_2}{x(x+N+d)(2x-1+M+d)_2},\\
  \text{Family (iii)}:\ \ \bar{B}(x)
  &=B(x+M)\,\frac{q^M(q^{x-N}\,;q)_M}{(q^{x+1-N}\,;q)_M}=B(x+M)\,
  \frac{q^M(1-q^{x-N})}{1-q^{x+M-N}},\\
  \bar{D}(x)&=D(x)\,\frac{q^{-M}(q^{x+1}\,;q)_M}{(q^{x}\,;q)_M}
  =D(x)\,\frac{q^{-M}(1-q^{x+M})}{1-q^x},\\
  \text{Family (iv)}:\ \ \bar{B}(x)
  &=B(x+M)\,\frac{(q^{x-N}\,;q)_M}{(q^{x+1-N}\,;q)_M}=B(x+M)\,
  \frac{1-q^{x-N}}{1-q^{x+M-N}},\\
  \bar{D}(x)&=D(x)\,\frac{(q^{x+1}\,;q)_M}{(q^{x}\,;q)_M}
  =D(x)\,\frac{1-q^{x+M}}{1-q^x},\\
  \text{Family (v)}:\ \ \bar{B}(x)
  &=B(x+M)\,\frac{(1-q^{x-N})(1-dq^x)(dq^{2x+2M}\,;q)_2}
  {(1-q^{x+M-N})(1-dq^{x+M})(dq^{2x+M}\,;q)_2},\\
  \bar{D}(x)&=D(x)\,\frac{(1-q^{x+M})(1-dq^{x+M+N})(dq^{2x-1}\,;q)_2}
  {(1-q^x)(1-dq^{x+N})(dq^{2x+M-1}\,;q)_2}.
\end{align}
It is straightforward to verify {\bf Theorem\,\ref{theo:sim}} for the
transformations of $B(x)$ and $D(x)$ for the families (i)--(v) by consulting
the explicit expressions of $B(x)$ and $D(x)$ of each polynomial listed in
\S\,\ref{sec:data}. The transformations rules of the polynomials \eqref{f14P},
\eqref{f2P}, \eqref{f3P}, and \eqref{f5P} are the direct consequences of
those for $\bar{B}$ and $\bar{D}$.
We list the explicit forms of the transformations of the polynomials of
degree $n\in\cX$ in the following:
\begin{theo}
\label{theo:42}
\begin{align}
  \text{Family {\rm (i)}}:\ \ &
  \sum_{j=0}^M(-1)^j\binom{M}{j}(x+1+j)_{M-j}(x-N)_j
  \check{P}_n(x+j;N,{\bm{\lambda}})\n
  &=(N+1)_M\,\check{P}_n(x+M;N+M,{\bm{\lambda}}),
  \label{f1eqF}\\
  \text{Family {\rm (ii)}}:\ \ &\sum_{j=0}^M(-1)^j\binom{M}{j}
  \mathcal{T}(x,M,j,d)(x+1+j,x+1+j+N+d)_{M-j}\n
  &\quad\times(x-N,x+d)_j\,\check{P}_n(x+j;N,d,\bar{\bm{\lambda}})\n
  &=(N+1)_M(2x+1+d)_{2M-1}\,\check{P}_n(x+M;N+M,d-M,\bar{\bm{\lambda}}),
  \label{f2eqF}\\
  \text{Family {\rm (iii)}}:\ \ &
  \sum_{j=0}^M(-1)^j\genfrac{[}{]}{0pt}{}{\,M\,}{j}q^{\frac12j(j+1)+M(N-j)}
  (q^{x+1+j}\,;q)_{M-j}(q^{x-N}\,;q)_j\check{P}_n(x+j;N,p)\n
  &=(q^{N+1}\,;q)_Mq^{Mx}\,\check{P}_n(x+M;N+M,pq^{-M}),
  \label{f3eqF}\\
  \text{Family {\rm (iv)}}:\ \ &
  \sum_{j=0}^M(-1)^j\genfrac{[}{]}{0pt}{}{\,M\,}{j}
  q^{\frac12j(j+1)+Nj}(q^{x+1+j}\,;q)_{M-j}(q^{x-N}\,;q)_j
  \check{P}_n(x+j;N,{\bm{\lambda}})\n
  &=(q^{N+1}\,;q)_M\,\check{P}_n(x+M;N+M,{\bm{\lambda}}),
  \label{f4eqF}\\
  \text{Family {\rm (v)}}:\ \ &
  \sum_{j=0}^M(-1)^j\genfrac{[}{]}{0pt}{}{\,M\,}{j}
  q^{\frac12j(j+1)+Nj}\mathcal{T}_q(x,M,j,d)
  (q^{x+1+j},dq^{x+1+j+N}\,;q)_{M-j}\n
  &\quad\times(q^{x-N},dq^x\,;q)_j\,\check{P}_n(x+j;N,d,\bar{\bm{\lambda}})\n
  &=(q^{N+1}\,;q)_M(dq^{2x+1}\,;q)_{2M-1}
  \,\check{P}_n(x+M;N+M,dq^{-M},\bar{\bm{\lambda}}).
  \label{f5eqF}
\end{align}
\end{theo}
The summations in l.h.s.\ are taken from \eqref{1eqform}, \eqref{2eqform},
\eqref{3eqform}, \eqref{4eqform} and \eqref{5eqform}. The factors in r.h.s.\ are
determined by setting $n=0$. A consistency check of the formulas can be
done by setting $x=-M$.
Again, it is straightforward to verify these transformation formulas for the
families (i)--(v) by consulting the explicit expressions of the polynomials.
Since the size parameter of the matrix is now $N+M$, there are additional
members of the orthogonal polynomials. They are simply
$\{\check{P}_{N+1+m}(x+M,N+M,{\bm{\lambda}})\}$ $(0\leq m\leq M)$, for
family (i), etc. This concludes the proof of {\bf Theorem\,\ref{theo:sim}}.

\bigskip
It is interesting and instructive to scrutinise {\bf Theorem\,\ref{theo:42}}
from a different perspective. The formulas \eqref{f1eqF}--\eqref{f5eqF} contain
$\sum_{j=0}^M$, $\binom{M}{j}$ or $\genfrac{[}{]}{0pt}{}{\,M\,}{j}$ and
subscripts $M-j$ and $j$, showing a semblance to binomial expansions.
Let us introduce operators $\tilde{\mathcal{F}}$ to realise one step
transformation for each family.
For $M=1$, {\bf Theorem\,\ref{theo:42}} tells
\begin{align}
  \text{(i)}:
  &\ \ \cFt(x,N)\eqdef\frac{1}{N+1}\bigl(x+1+(N-x)e^{\partial}\bigr),\\
  &\qquad\cFt(x,N)\check{P}_n(x;N,{\bm{\lambda}})
  =\check{P}_n(x+1;N+1,{\bm{\lambda}}),\n
  \text{(ii)}:
  &\ \ \cFt(x,N,d)\eqdef\frac{1}{(N+1)(2x+1+d)}
  \bigl((x+1)(x+1+N+d)+(N-x)(x+d)e^{\partial}\bigr),\\
  &\qquad\cFt(x,N,d)\check{P}_n(x;N,d,\bar{\bm{\lambda}})
  =\check{P}_n(x+1;N+1,d-1,\bar{\bm{\lambda}}),\n 
  \text{(iii)}:
  &\ \ \cFt(x,N)\eqdef\frac{q^{N-x}}{1-q^{N+1}}
  \bigl(1-q^{x+1}+(q^{x-N}-1)e^{\partial}\bigr),\\
  &\qquad\cFt(x,N)\check{P}_n(x;N,p)=\check{P}_n(x+1;N+1,pq^{-1}),\n 
  \text{(iv)}:
  &\ \ \cFt(x,N)\eqdef\frac{1}{1-q^{N+1}}
  \bigl(1-q^{x+1}+q^{N+1}(q^{x-N}-1)e^{\partial}\bigr),\\
  &\qquad\cFt(x,N)\check{P}_n(x;N,{\bm{\lambda}})
  =\check{P}_n(x+1;N+1,{\bm{\lambda}}),\n
  \text{(v)}:
  &\ \ \cFt(x,N,d)\eqdef\frac{1}{(1-q^{N+1})(1-dq^{2x+1})}\n
  &\phantom{\ \ \cFt(x,N,d)\eqdef}\times
  \bigl((1-q^{x+1})(1-dq^{x+1+N})+q^{N+1}(q^{x-N}-1)(1-dq^x)e^{\partial}\bigr),\\
  &\qquad\cFt(x,N,d)\check{P}_n(x;N,d,\bar{\bm{\lambda}})
  =\check{P}_n(x+1;N+1,dq^{-1},\bar{\bm{\lambda}}). 
  \nonumber
\end{align}
Here $e^\partial$ is the finite shift operator \eqref{finshift}.
Let us tentatively call $\cFt$ a forward $x$-shift operator.
Each $\cFt$ has two terms, one is an ordinary function and the other is a
function times $e^{\partial}$. Since the polynomial is mapped from
$\check{P}_n(x;N,{\bm{\lambda}})$ to $\check{P}_n(x+1;N+1,{\bm{\lambda}})$,
for family (i), the next steps are obviously
\begin{align*}
  &\text{(i)}:\ \ \cFt(x+1,N+1)\cFt(x,N)\check{P}_n(x;N,{\bm{\lambda}})
  =\check{P}_n(x+2;N+2,{\bm{\lambda}}),\\
  &\phantom{\text{(i)}:\ \ }\cFt(x+2,N+2)\cFt(x+1,N+1)\cFt(x,N)
  \check{P}_n(x;N,{\bm{\lambda}})=\check{P}_n(x+3;N+3,\bm{\lambda}),\\
  &\hspace{4cm}\vdots\ ,
\end{align*}
leading to
\begin{alignat}{2}
  \text{(i),\,(iv)}:
  &&\ \ \laprod{k=0}{M-1}\cFt(x+k,N+k)\cdot\check{P}_n(x;N,\bm{\lambda})
  &=\check{P}_n(x+M;N+M,\bm{\lambda}),\\
  \text{(ii)}:
  &&\ \laprod{k=0}{M-1}\cFt(x+k,N+k,d-k)\cdot
  \check{P}_n(x;N,d,\bar{\bm{\lambda}})
  &=\check{P}_n(x+M;N+M,d-M,\bar{\bm{\lambda}}),\!\\
  \text{(iii)}:
  &&\ \ \laprod{k=0}{M-1}\cFt(x+k,N+k)\cdot\check{P}_n(x;N,p)
  &=\check{P}_n(x+M;N+M,pq^{-M}),\\
  \text{(v)}:
  &&\ \ \laprod{k=0}{M-1}\cFt(x+k,N+k,dq^{-k})\cdot
  \check{P}_n(x;N,d,\bar{\bm{\lambda}})
  &=\check{P}_n(x+M;N+M,dq^{-M},\bar{\bm{\lambda}}).
\end{alignat}
Here we use the convention to express the ordered product
\begin{equation}
  \laprod{j=1}{n}a_j\eqdef a_n\cdots a_2a_1.
\end{equation}
Thus we arrive at the following
\begin{theo}
\label{theo:43}
The formulas representing the special cases of the multi-indexed polynomials
\eqref{f1eqF}--\eqref{f5eqF} are factorised,
\begin{align}
  {\rm (i)}:
  &\ \ \laprod{k=0}{M-1}\cFt(x+k,N+k)=\frac{1}{(N+1)_M}
  \sum_{j=0}^M(-1)^j\binom{M}{j}(x+1+j)_{M-j}(x-N)_j\,e^{j\partial},\\
  {\rm (ii)}:
  &\ \ \laprod{k=0}{M-1}\cFt(x+k,N+k,d-k)\n
  &=\frac{1}{(N+1)_M(2x+1+d)_{2M-1}}
  \sum_{j=0}^M(-1)^j\binom{M}{j}\mathcal{T}(x,M,j,d)\n
  &\qquad\qquad\times
  (x+1+j,x+1+j+N+d)_{M-j}(x-N,x+d)_j\,e^{j\partial},\\
  {\rm (iii)}:
  &\ \ \laprod{k=0}{M-1}\cFt(x+k,N+k)\n
  &=\frac{q^{-Mx}}{(q^{N+1}\,;q)_M}
  \sum_{j=0}^M(-1)^j\genfrac{[}{]}{0pt}{}{\,M\,}{j}q^{\frac12j(j+1)+M(N-j)}
  (q^{x+1+j}\,;q)_{M-j}(q^{x-N}\,;q)_j\,e^{j\partial},\\
  {\rm (iv)}:
  &\ \ \laprod{k=0}{M-1}\cFt(x+k,N+k)\n
  &=\frac{1}{(q^{N+1}\,;q)_M}
  \sum_{j=0}^M(-1)^j\genfrac{[}{]}{0pt}{}{\,M\,}{j}
  q^{\frac12j(j+1)+Nj}(q^{x+1+j}\,;q)_{M-j}(q^{x-N}\,;q)_j\,e^{j\partial},\\
  {\rm (v)}:
  &\ \ \laprod{k=0}{M-1}\cFt(x+k,N+k,dq^{-k})\n
  &=\frac{1}{(q^{N+1}\,;q)_M(dq^{2x+1}\,;q)_{2M-1}}
  \sum_{j=0}^M(-1)^j\genfrac{[}{]}{0pt}{}{\,M\,}{j}
  q^{\frac12j(j+1)+Nj}\mathcal{T}_q(x,M,j,d)\n
  &\qquad\qquad\times
  (q^{x+1+j},dq^{x+1+j+N}\,;q)_{M-j}(q^{x-N},dq^x\,;q)_j\,e^{j\partial}.
\end{align}
By expanding the l.h.s and collecting terms containing $e^{j\partial}$ gives
r.h.s. 
\end{theo}
These formulas can be demonstrated by straightforward induction.
For example, for (i),
\begin{align*}
  &\quad\laprod{k=0}M\cFt(x+k,N+k)
  =\cFt(x+M,N+M)\laprod{k=0}{M-1}\cFt(x+k,N+k)\\
  &=\frac{1}{(N+1)_{M+1}}
  \sum_{j=0}^M\Bigl((-1)^j\binom{M}{j}(x+1+j)_{M+1-j}(x-N)_j\,e^{j\partial}\\
  &\phantom{=\frac{1}{(N+1)_{M+1}}\sum_{j=0}^M\Bigl(}
  -(-1)^j\binom{M}{j}(x+2+j)_{M-j}(x-N)_{j+1}\,e^{(j+1)\partial}\Bigr)\\
  &=\frac{1}{(N+1)_{M+1}}
  \sum_{j=0}^{M+1}(-1)^j\binom{M+1}{j}
  (x+1+j)_{M+1-j}(x-N)_j\,e^{j\partial}.
\end{align*}
Here we use Pascal' triangle relation
\begin{equation}
  \binom{M}{j}+\binom{M}{j-1}=\binom{M+1}{j}.
  \label{binomid1}
\end{equation}
The $q$-analogues of the above relation
\begin{equation}
  \qbinom{M}{j}q^j+\qbinom{M}{j-1}=\qbinom{M+1}{j},\quad
  \qbinom{M}{j}+\qbinom{M}{j-1}q^{M+1-j}=\qbinom{M+1}{j}
  \label{qbinomid1}
\end{equation}
are used in the inductions for (iii)--(v) and the relations
\begin{equation}
  j\binom{M}{j}-(M+1-j)\binom{M}{j-1}=0,\quad
  (1-q^j)\qbinom{M}{j}-(1-q^{M+1-j})\qbinom{M}{j-1}=0
\end{equation}
are used for (ii) and (v).

It is well known that the difference operator $\widetilde{\mathcal{H}}$
\eqref{hteq}, \eqref{hteq1} is factorised into the forward $\mathcal{F}$ and
backward $\mathcal{B}$ shift operators \cite{kls,os12},
\begin{align}
  &\widetilde{\mathcal{H}}(N,\bm{\lambda})
  =B(x;N,\bm{\lambda})(1-e^\partial)+D(x;N,\bm{\lambda})(1-e^{-\partial})
  =\mathcal{B}(N,\bm{\lambda})\mathcal{F}(N,\bm{\lambda}),
  \label{BFfact}
 \end{align}
which is an expression of the shape invariance. For example, for the Racah, 
\begin{align}
  &\text{Racah}:\ \ \mathcal{F}(N,\bm{\lambda})
  =\frac{Nbc}{2x+d+1}(1-e^{\partial}),\n
  &\phantom{\text{Racah}:\ \ }\ \mathcal{B}(N,\bm{\lambda})
  =\frac1{Nbc}\bigl(B(x;N,\bm{\lambda})-D(x;N,\bm{\lambda})
  e^{-\partial}\bigr)(2x+d+1),\\
  &\mathcal{F}(N,\bm{\lambda})\check{P}_n(x;N,b,c,d)
  =\mathcal{E}(n;\bm{\lambda})\check{P}_{n-1}(x;N-1,b+1,c+1,d+1),\n
  &\mathcal{B}(N,\bm{\lambda})\check{P}_{n-1}(x;N-1,b+1,c+1,d+1)
  =\check{P}_n(x;N,b,c,d).
\end{align}
These $\mathcal{F}$ and $\mathcal{B}$ change the degree $n$ but keep $x$
unchanged, in contrast with $\tilde{\mathcal{F}}$, which changes $x$ and keeps
the degree $n$ unchanged.
In a similar fashion as above \eqref{BFfact}, reflecting the new type of shape
invariance, $\widetilde{\mathcal{H}}-\mathcal{E}(N+1)$ is factorised by the
forward and backward $x$-shift operators,
\begin{equation}
  \widetilde{\mathcal{H}}(N,\bm{\lambda})-\mathcal{E}(N+1;\bm{\lambda})
  =-\cBt(x,N,\bm{\lambda})\cFt(x,N,\bm{\lambda}),
  \label{tH2}
\end{equation}
as the zero norm solutions start at degree $N+1$.
The explicit forms of the operator $\tilde{\mathcal{B}}$ are
\begin{alignat}{2}
  \text{K}:
  &\ \ &\cBt(x,N,p)&\eqdef(N+1)\bigl(p+(1-p)e^{-\partial}\bigr),\\
  \text{H}:
  &\ \ &\cBt(x,N,\bm{\lambda})&\eqdef(N+1)\bigl(x+a+(b+N-x)e^{-\partial}\bigr),\\
  \text{R}:
  &\ \ &\cBt(x,N,\bm{\lambda})&\eqdef(N+1)\Bigl(\frac{(x+b)(x+c)}{2x+d}
  +\frac{(b-d-x)(x+d-c)}{2x+d}e^{-\partial}\Bigr),\\
  \text{dH}:
  &\ \ &\cBt(x,N,\bm{\lambda})&\eqdef(N+1)\Bigl(\frac{x+a}{2x-1+a+b}
  +\frac{x+b-1}{2x-1+a+b}e^{-\partial}\Bigr),\\
  \text{dq$q$K}:
  &\ \ &\cBt(x,N,p)&\eqdef(1-q^{N+1})\bigl(p^{-1}q^{-x-N-1}
  +q^{-N-1}(1-p^{-1}q^{-x})e^{-\partial}\bigr),\\
  \text{$q$H}:
  &\ \ &\cBt(x,N,\bm{\lambda})&\eqdef(1-q^{N+1})\bigl(q^{-N-1}(1-aq^x)
  +aq^{-1}(q^{x-N}-b)e^{-\partial}\bigr),\\
  \text{$q$K}:
  &\ \ &\cBt(x,N,p)&\eqdef(1-q^{N+1})\bigl(q^{-N-1}
  +p\,e^{-\partial}\bigr),\\
  \text{q$q$K}:
  &\ \ &\cBt(x,N,p)&\eqdef(1-q^{N+1})\bigl(p^{-1}q^{x-N-1}
  +(1-p^{-1}q^{x-N-1})e^{-\partial}\bigr),\\
  \text{a$q$K}:
  &\ \ &\cBt(x,N,p)&\eqdef(1-q^{N+1})\bigl(q^{-N-1}(1-pq^{x+1})
  +pq^{x-N}e^{-\partial}\bigr),\\
  \text{$q$R}:
  &\ \ &\cBt(x,N,\bm{\lambda})&\eqdef(1-q^{N+1})\Bigl(
  \frac{q^{-N-1}(1-bq^x)(1-cq^x)}{1-dq^{2x}}\n
  &&&\phantom{\eqdef(1-q^{N+1})\Bigl(}
  +\tilde{d}\,\frac{(b^{-1}dq^x-1)(1-c^{-1}dq^x)}{1-dq^{2x}}
  e^{-\partial}\Bigr),\\
  \text{d$q$H}:
  &\ \ &\cBt(x,N,\bm{\lambda})&\eqdef(1-q^{N+1})\Bigl(
  \frac{q^{-N-1}(1-aq^x)}{1-abq^{2x-1}}
  +\frac{aq^{x-N-1}(1-bq^{x-1})}{1-abq^{2x-1}}e^{-\partial}\Bigr),\\
  \text{d$q$K}:
  &\ \ &\cBt(x,N,p)&\eqdef(1-q^{N+1})\Bigl(
  \frac{q^{-N-1}}{1+pq^{2x}}
  +\frac{pq^{2x-N-1}}{1+pq^{2x}}e^{-\partial}\Bigr).
\end{alignat}
The action of these operators on the polynomials is as follows:
\begin{alignat}{2}
  &\text{(i),\,(iv)}:
  &\ \ \cBt(x,N)\check{P}_n(x+1;N+1,\bm{\lambda})
  &=\bigl(\mathcal{E}(N+1;\bm{\lambda})-\mathcal{E}(n;\bm{\lambda})\bigr)
  \check{P}_n(x;N,\bm{\lambda}),\\
  &\text{(ii)}:
  &\hspace*{-6mm}\cBt(x,N,d,\bar{\bm{\lambda}})
  \check{P}_n(x+1;N+1,d-1,\bar{\bm{\lambda}})
  &=\bigl(\mathcal{E}(N+1;\bm{\lambda})-\mathcal{E}(n;\bm{\lambda})\bigr)
  \check{P}_n(x;N,d,\bar{\bm{\lambda}}),\!\!\!\\
  &\text{(iii)}:
  &\ \ \cBt(x,N,p)\check{P}_n(x+1;N+1,pq^{-1})
  &=\bigl(\mathcal{E}(N+1)-\mathcal{E}(n)\bigr)
  \check{P}_n(x;N,p),\\
  &\text{(v)}:
  &\hspace*{-6mm}\cBt(x,N,d,\bar{\bm{\lambda}})
  \check{P}_n(x+1;N+1,dq^{-1},\bar{\bm{\lambda}})
  &=\bigl(\mathcal{E}(N+1;\bm{\lambda})-\mathcal{E}(n;\bm{\lambda})\bigr)
   \check{P}_n(x;N,d,\bar{\bm{\lambda}}).\!\!\!
\end{alignat}

\section{Summary and Comments}
\label{sec:comm}

The reported ``Diophantine'' and factorisation properties of the Jacobi,
Laguerre etc \cite{cal}, Wilson and Askey-Wilson polynomials etc
\cite{chen-ismail} are very interesting but rather puzzling, as these
properties require the parameter ranges in which the orthogonality does not
hold.
In this paper we show, after the seminal work of \cite{miki}, that these
properties are shared by all the finite polynomials in the Askey scheme in
the conventional orthogonality parameter ranges. All the {\em monic} higher
degree polynomials $\{\check{P}_{N+1+m}(x)\}$ ($m\in\mathbb{Z}_{\ge0}$) in
this group show the ``Diophantine'' and factorisation properties
{\bf Theorem\,\ref{theo1}} \eqref{Dfac}, \eqref{lamdef}. Here $N$ is the
maximal degree of the corresponding non-monic orthogonal polynomials.
A simple and intuitive explanation is that these higher degree monic
polynomials are the zero-norm solution of certain tri-diagonal real symmetric
matrix $\mathcal{H}$ \eqref{tdefcomp}. The explicit expressions of the
``Diophantine'' and factorisation for the twelve polynomials are presented in
\S\,\ref{sec:zero}. These higher degree monic polynomials can be used as seed
solutions for generating new types multi-indexed polynomials based on the
twelve finite orthogonal polynomials. In order to obtain genuine multi-indexed
orthogonal polynomials, correct choices of the seed solutions and the
parameter ranges are essential. A simplest choice of $M$ contiguous lowest
degree $\{\check{P}_{N+1+m}(x)\}$ ($m=0,1,\ldots,M-1$) generate proper
multi-indexed orthogonal polynomials as shown in \S\,\ref{sec:shape} for each
of the twelve polynomials listed in \S\,\ref{sec:data}.

In paper $\I$\,\cite{cal} Calogero and coauthors reported ``Diophantine'' and
factorisation properties of some finite polynomials, the Racah, Hahn, dual
Hahn and Krawtchouk. In these cases, the phenomena occur only at the
particular degree with which the parameters are tuned.
In `Remarks' of paper $\II$\,\cite{chen-ismail}, Chen and Ismail mentioned
a possible explanation of the ``Diophantine'' and factorisation properties
of the Wilson polynomial as the occurrence of certain discrete masses in the
orthogonality measure when some of the parameters are negative.
It is a good challenge to provide a perspective in which their explanation
and ours could be unified, since the Wilson and Askey-Wilson polynomials are
the most general ones in the Askey scheme.

A few remarks on the results of Miki-Tsujimoto-Vinet paper \cite{miki}.
They reported four types of single-indexed exceptional Krawtchouk polynomials.
Among the them, the first one, using the polynomial itself as the seed
solution, is the well-known one called Krein-Adler polynomials, as briefly
mentioned in \S\,\ref{sec:set}. The second type is the $M=1$ case of the
multi-indexed polynomials discussed in \S\,\ref{sec:newmulti}. The remaining
two types exist only for the Krawtchouk polynomials. They are obtained from
the first and second types by using the mirror symmetry \eqref{Kpmirror} of
the Krawtchouk. Another polynomial in the same $\eta(x)=x$ family, the Hahn
polynomial, has similar mirror symmetry \eqref{Hpd}. The forms of the other
$\eta(x)$, family (ii)--(v), are not compatible with the mirror symmetry.

\section*{Acknowledgements}

S.\,O. is supported by JSPS KAKENHI Grant Number JP19K03667.



\begin{thebibliography}{99}

\bibitem{cal}
M.\,Bruschi, F.\,Calogero and R.\,Droghei,
``Additional recursion relations, factorizations, and Diophantine properties
associated with the polynomials of the Askey scheme,''
Adv. Math. Phys. 2009 (2009) 268134 (43pp).
 
\bibitem{chen-ismail}
Y.\,Chen and M.\,Ismail,
``Hypergeometric origins of Diophantine properties associated with the Askey
scheme,''
Proceedings of the American Mathematical Society {\bf 138} (2010) 943-951.

\bibitem{askey}
G.\,E.\,Andrews, R.\,Askey and R.\,Roy,
{\it Special Functions},
Encyclopedia of mathematics and its applications,
Cambridge Univ. Press, Cambridge, (1999).

\bibitem{ismail}
M.\,E.\,H.\,Ismail,
{\it Classical and Quantum Orthogonal Polynomials in One Variable\/},
Encyclopedia of mathematics and its applications,
Cambridge Univ. Press, Cambridge, (2005).

\bibitem{kls}
R.\,Koekoek, P.\,A.\,Lesky and R.\,F.\,Swarttouw,
{\it Hypergeometric orthogonal polynomials and their $q$-analogues,\/}
Springer Monographs in Mathematics,
Springer-Verlag Berlin-Heidelberg, (2010).

\bibitem{gasper}
G.\,Gasper and M.\,Rahman,
{\it Basic hypergeometric series\/}, 2nd ed.,
Encyclopedia of mathematics and its applications,
Cambridge Univ. Press, Cambridge, (2004).

\bibitem{os12}
S.\,Odake and R.\,Sasaki,
``Orthogonal Polynomials from Hermitian Matrices,''
J. Math. Phys. {\bf 49} (2008) 053503 (43 pp),
{\tt arXiv:0712.4106[math.CA]}.

\bibitem{nikiforov}
A.F.\,Nikiforov, S.K.\,Suslov and V.B.\,Uvarov,
{\it Classical Orthogonal Polynomials of a Discrete Variable\/},
Springer-Verlag, Berlin, (1991).

\bibitem{miki}
H.\,Miki,\ S.\,Tsujimoto and L.\,Vinet,
``The single-indexed exceptional Krawtchouk polynomials,''
{\tt arXiv:2201.12359[math.CA]}.

\bibitem{crum}
M.\,M.\,Crum,
``Associated Sturm-Liouville systems,''
Quart. J. Math. Oxford Ser. (2) {\bf 6} (1955) 121-127,
{\tt arXiv:physics/9908019}.

\bibitem{os22}
S.\,Odake and R.\,Sasaki,
``Dual Christoffel transformations,''
Prog. Theor. Phys. {\bf 126} (2011) 1-34.
{\tt arXiv:1101.5468[math-ph]}.

\bibitem{krein}
M.\,G.\,Krein,
Doklady Acad. Nauk. CCCP, {\bf 113} (1957) 970-973.

\bibitem{adler}
V.\,\'E.\,Adler,
``A modification of Crum's method,''
Theor. Math. Phys. {\bf 101} (1994) 1381-1386.

\bibitem{os26}
S.\,Odake and R.\,Sasaki,
``Multi-indexed ($q$-)Racah polynomials,''
J. Phys. {\bf A45} (2012) 385201 (21 pp),
{\tt arXiv:1203.5868[math-ph]}.

\bibitem{gomez}
D.\,Gomez-Ullate, N.\,Kamran and R.\,Milson,
``An extended class of orthogonal polynomials defined by a Sturm-Liouville
problem,''
J. Math. Anal. Appl. {\bf 359} (2009) 352-367, {\tt arXiv:0807.3939[math-ph]};
``An extension of Bochner's problem: exceptional invariant sub-spaces,''
J. Approx. Theory {\bf 162} (2010) 987-1006, {\tt arXiv:0805.3376[math-ph]}.

\bibitem{quesne}
C.\,Quesne,
``Exceptional orthogonal polynomials, exactly solvable potentials and
supersymmetry,''
J. Phys. A: Math. Theor. {\bf 41} (2008) 392001 (6 pp),
{\tt arXiv:0807.4087\hspace{0pt}[quant-ph]}.

\bibitem{os16}
S.\,Odake and R.\,Sasaki,
``Infinitely many shape invariant potentials and new orthogonal polynomials,''
Phys. Lett. {\bf B679} (2009) 414-417,
{\tt arXiv:0906.0142[math-ph]}.

\bibitem{os25}
S.\,Odake and R.\,Sasaki,
``Exactly solvable quantum mechanics and infinite families of multi-indexed
orthogonal polynomials,''
Phys. Lett. {\bf B702} (2011) 164-170,
{\tt arXiv:1105.\hspace{0pt}0508[math-ph]}.

\bibitem{os24}
S.\,Odake and R.\,Sasaki,
``Discrete quantum mechanics,'' (Topical Review)
J. Phys. {\bf A44} (2011) 353001 (47 pp),
{\tt arXiv:1104.0473[math-ph]}.

\bibitem{genden}
L.\,E.\,Gendenshtein,
``Derivation of exact spectra of the Schroedinger equation by means of
supersymmetry,''
JETP Lett. {\bf 38} (1983) 356-359.

\end{thebibliography}
\end{document}